\documentclass[11pt]{amsart}%%%%%,leqno
\usepackage{amsmath,amssymb, graphics, amscd,latexsym}

\makeatletter

\newtheorem{Theorem}{Theorem}%[section]
\newtheorem{Lemma}[Theorem]{Lemma}
\newtheorem{Corollary}[Theorem]{Corollary}
\newtheorem{Proposition}[Theorem]{Proposition}

\theoremstyle{remark}
\newtheorem{Remark}[Theorem]{Remark}

\newtheorem{Observation}[Theorem]{Observation}

\begin{document}
\newcommand{\eps}{\varepsilon}
\newcommand{\om}{\omega}
\newcommand\Om{\Omega}
\newcommand\la{\lambda}
\newcommand\vphi{\varphi}
\newcommand\vrho{\varrho}
\newcommand\al{\alpha}
\newcommand\La{\Lambda}
\newcommand\si{\sigma}
\newcommand\be{\beta}
\newcommand\Si{\Sigma}
\newcommand\ga{\gamma}
\newcommand\Ga{\Gamma}
\newcommand\de{\delta}
\newcommand\De{\Delta}

\newcommand\cA{\mathcal  A}
\newcommand\cB{\mathcal B}
\newcommand\cD{\mathcal  D}
\newcommand\cM{\mathcal  M}
\newcommand\cN{\mathcal  N}
\newcommand\cT{\mathcal  T}
\newcommand\cP{\mathcal  P}
\newcommand\cp{\mathcal  p}
\newcommand\cQ{\mathcal  Q}
\newcommand\cG{\mathcal G}
\newcommand\cq{\mathcal  q}
\newcommand\cc{\mathcal  c}
\newcommand\cs{\mathcal  s}
\newcommand\cS{\mathcal  S}
\newcommand\ct{\mathcal  t}
\newcommand\cZ{\mathcal  Z}
\newcommand\cR{\mathcal  R}
\newcommand\cu{\mathcal  u}
\newcommand\cU{\mathcal  U}
\newcommand\cI{\mathcal  I}
\newcommand\cJ{\mathcal  J}
\newcommand\co{\mathcal  o}
\newcommand\cO{\mathcal  O}
\newcommand\cv{\mathcal  v}
\newcommand\cV{\mathcal  V}
\newcommand\cx{\mathcal  x}
\newcommand\cX{\mathcal  X}
\newcommand\cw{\mathcal  w}
\newcommand\ck{\mathcal  k}
\newcommand\cK{\mathcal  K}
\newcommand\cW{\mathcal  W}
\newcommand\cz{\mathcal  z}
\newcommand\cy{\mathcal  y}
\newcommand\ca{\mathcal  a}
\newcommand\ch{\mathcal  h}
\newcommand\cH{\mathcal  H}
\newcommand\cF{\mathcal F}
\newcommand\bfG{\mbox {\bf  G}}
\newcommand\bfg{\mbox {\bf  g}}
\newcommand\bfC{\mbox {\bf  C}}
\newcommand\bfN{\mbox {\bf  N}}
\newcommand\bfT{\mbox {\bf  T}}
\newcommand\bfP{\mbox {\bf  P}}
\newcommand\bfp{\mbox {\bf  p}}
\newcommand\bfQ{\mbox {\bf  Q}}
\newcommand\bfq{\mbox {\bf  q}}
\newcommand\bfc{\mbox {\bf  c}}
\newcommand\bfs{\mbox {\bf  s}}
\newcommand\bfS{\mbox {\bf  S}}
\newcommand\bft{\mbox {\bf  t}}
\newcommand\bfZ{\mbox {\bf  Z}}
\newcommand\bfR{\mbox {\bf  R}}
\newcommand\bfu{\mbox {\bf  u}}
\newcommand\bfU{\mbox {\bf  U}}
\newcommand\bfo{\mbox {\bf  o}}
\newcommand\bfO{\mbox {\bf  O}}
\newcommand\bfv{\mbox {\bf  v}}
\newcommand\bfV{\mbox {\bf  V}}
\newcommand\bfx{\mbox {\bf  x}}
\newcommand\bfX{\mbox {\bf  X}}
\newcommand\bfw{\mbox {\bf  w}}
\newcommand\bfk{\mbox {\bf  k}}
\newcommand\bfK{\mbox {\bf  K}}
\newcommand\bfW{\mbox {\bf  W}}
\newcommand\bfz{\mbox {\bf  z}}
\newcommand\bfy{\mbox {\bf  y}}
\newcommand\bfa{\mbox {\bf  a}}
\newcommand\bfh{\mbox {\bf  h}}
\newcommand\bfH{\mbox {\bf  H}}
\newcommand\bfJ{\mbox {\bf  J}}
\newcommand\bfj{\mbox {\bf  j}}
\newcommand\bbC{\mbox {\mathbb C}}
\newcommand\bbN{\mbox {\mathbb N}}
\newcommand\bbT{\mbox {\mathbb T}}
\newcommand\bbP{\mbox {\mathbb P}}
\newcommand\bbQ{\mbox {\mathbb Q}}
\newcommand\bbS{\mbox {\mathbb S}}
\newcommand\bbZ{\mbox {\mathbb Z}}
\newcommand\bbR{\mbox {\mathbb R}}
\newcommand\bbU{\mbox {\mathbb U}}
\newcommand\bbO{\mbox {\mathbb O}}
\newcommand\bbV{\mbox {\mathbb V}}
\newcommand\bbX{\mbox {\mathbb X}}
\newcommand\bbK{\mbox {\mathbb K}}
\newcommand\bbW{\mbox {\mathbb W}}
\newcommand\bbH{\mbox {\mathbb H}}

\newcommand\apeq{\fallingdotseq}
\newcommand\Lrarrow{\Leftrightarrow}
\newcommand\bij{\leftrightarrow}
\newcommand\Rarrow{\Rightarrow}
\newcommand\Larrow{\Leftarrow}
\newcommand\nin{\noindent}
\newcommand\ninpar{\par \noindent}
\newcommand\nlind{\nl \indent}
\newcommand\nl{\newline}
\newcommand\what{\widehat}
\newcommand\tl{\tilde}
\newcommand\wtl{\widetilde}
\newcommand\order{\mbox{\text{order}\/}}
\newcommand\GL{\text{GL}\/}
\newcommand\PGL{\text{PGL}\/}
\newcommand\Spec{\text{Spec}\/}
\newcommand\weight{\text{weight}\/}
\newcommand\ord{\text{ord}\/}
\newcommand\Int{\text{Int}\/}
\newcommand\grad{\text{grad}\/}
\newcommand\Ind{\text{Ind}\/}
\newcommand\Disc{\text{Disc}\/}
\newcommand\Ker{\text{Ker}\/}
\newcommand\Image{\text{Image}\/}
\newcommand\Coker{\text{Coker}\/}
\newcommand\Id{\text{Id}\/}
\newcommand\id{\text{id}}
\newcommand\dsum{\text{\amalg}}
\newcommand\val{\text{val}}
\newcommand\mv{\text{m-vector}}
\newcommand\iv{\text{i-vector}}
\newcommand\minimum{\text{minimum}\/}
\newcommand\modulo{\text{modulo}\/}
\newcommand\Aut{\text{Aut}\/}
\newcommand\PSL{\text{PSL}\/}
\newcommand\Res{\text{Res}\/}
\newcommand\rank{\text{rank}\/}
\newcommand\codim{\text{codim}\/}
\newcommand\msdim{\text{ms-dim}\/}
\newcommand\icodim{\text{i-codim}\/}
\newcommand\ocodim{\text{o-codim}\/}
\newcommand\emdim{\text{ems-dim}\/}

\newcommand\Cone{\text{Cone}\/}
\newcommand\maximum{\text{maximum}\/}
\newcommand\Vol{\text{Vol}\/}
\newcommand\Coeff{\text{Coeff}\/}
\newcommand\lcm{\text{lcm}\/}
\newcommand\degree{\text{degree}\/}
\newcommand\Pol{\cal {POL}}
\newcommand\Ts{Tschirnhausen}
\newcommand\TS{Tschirnhausen approximate}
\newcommand\Stab{\text{Stab}\/}
\newcommand\civ{complete intersection variety}
\newcommand\nipar{\par \noindent}
\newcommand\wsim{\overset{w}{\sim}}
\newcommand\Gr{\bfZ_2*\bfZ_3}
\newcommand\QED{~~Q.E.D.}
\newcommand\bsq{$\blacksquare$}
\newcommand\bff{\mbox {\bf  f}}
\newcommand\newcommandby{:=}
\newcommand\inv{^{-1}}
\newcommand\nnt{(\text{nn-terms})}
\renewcommand{\subjclassname}{\textup{2000} Mathematics Subject Classification}
\def\mapright#1{\smash{\mathop{\longrightarrow}\limits^{{#1}}}}
%\def\maprightt#1#2{\smash{\mathop{\longrightarrow}\limits^{#1}}}
%左矢印
\def\mapleft#1{\smash{\mathop{\longleftarrow}\limits^{{#1}}}}
%同下向き矢印
\def\mapdown#1{\Big\downarrow\rlap{$\vcenter{\hbox{$#1$}}$}}
\def\mapdownn#1#2{\llap{$\vcenter{\hbox{$#1$}}$}\Big\downarrow\rlap{$\vcenter{\hbox{$#2$}}$}}
%同上向き矢印
\def\mapup#1{\Big\uparrow\rlap{$\vcenter{\hbox{$#1$}}$}}
\def\mapupp#1#2{\llap{$\vcenter{\hbox{$#1$}}$}\Big\uparrow\rlap{$\vcenter{\hbox{$#2$}}$}}

%同右下向き
\def\rdown#1{\searrow\rlap{$\vcenter{\hbox{$#1$}}$}}
\def\semap#1{\searrow\rlap{$\vcenter{\hbox{$#1$}}$}}

%同右上向き
\def\rup#1{\nearrow\rlap{$\vcenter{\hbox{$#1$}}$}}
\def\nemap#1{\nearrow\rlap{$\vcenter{\hbox{$#1$}}$}}
%同左下向き
\def\ldown#1{\swarrow\rlap{$\vcenter{\hbox{$#1$}}$}}
\def\swmap#1{\swarrow\rlap{$\vcenter{\hbox{$#1$}}$}}
%同左上向き
\def\lup#1{\nwarrow\rlap{$\vcenter{\hbox{$#1$}}$}}
\def\nwmap#1{\nwarrow\rlap{$\vcenter{\hbox{$#1$}}$}}
\def\defby{:=}
\def\eqby#1{\overset {#1}\to =}
\def\inv{^{-1}}
\def\bnu{{(\nu)}}

\title[Classification of sextics of  torus type
%{\it Draft: \today}
]
{Classification of sextics of  torus type}

\author
%Normally smooth divisors  {\it Draft: \today}
[M. Oka {\tiny and} D.T. Pho]
{Mutsuo Oka {\tiny and} Duc Tai Pho}

\address{\vtop{
\hbox{Mutsuo Oka}
\hbox{Department of Mathematics}
\hbox{Tokyo Metropolitan University}
\hbox{1-1 Mimami-Ohsawa, Hachioji-shi}
\hbox{Tokyo 192-0397}
\hbox{\rm{E-mail}:{\rm oka@comp.metro-u.ac.jp}}}
\hspace{1cm}
\vtop{\hbox{Duc Tai Pho}
\hbox{Dept. of Math. and Computer Sci.}
\hbox{Bar-Ilan University}
\hbox{52900 Ramat-Gan, Israel}
\hbox{\rm phoduct@macs.biu.ac.il}}}

%\thanks{}
\keywords{Torus type, maximal curves, dual curves}
\subjclass{14H10,14H45, 32S05.}
\maketitle

\begin{center}{Dedicated to Professor Tatsuo Suwa on his 60th
 birthday}
\end{center}

\begin{abstract}
In \cite{Pho}, the second author classified  configurations  of the
 singularities on
  tame sextics of torus type.
In this paper, we give a complete
classification of the singularities on irreducible sextics of  torus
 type,
without assuming the tameness of the sextics.
We show that there exists 121  configurations and
 there are 5 pairs and a triple of configurations
for which the corresponding moduli spaces coincide, ignoring the
 respective
torus decomposition.
\end{abstract}

\pagestyle{headings}
\section{Introduction}
We consider an irreducible sextics of torus type $C$ defined by 
\begin{eqnarray}\label{fixed-torus-decomp}
 C:\{(X,Y,Z)\in \bfP^2; F_2(X,Y,Z)^3+F_3(X,Y,Z)^2=0\}
\end{eqnarray}
where $F_i(X,Y,Z)$ is a homogeneous polynomial of degree $i$, for
$i=2,3$.  We consider the conic $C_2=\{F_2(X,Y,Z)=0 \}$ and 
the cubic $C_3=\{F_3(X,Y,Z)=0\}$.
% defined by $F_i(X,Y,Z)=0$ respectively.
 Let $\Si(C)$ be the set of singular points
of $C$. A singular point $P\in \Si(C)$ is called
{\em inner} (respectively {\em outer})
with respect to the given torus decomposition
(\ref{fixed-torus-decomp}) if $P\in C_2$ (resp. $P\notin C_2$).
We say that  $C$ is {\em tame} if $\Si(C)\subset C_2\cap C_3$.
For  tame sextics of  torus type, there are 25 
local  singularity types among which 20 appear on irreducible 
sextics of torus type by \cite{Pho}.
As   global singularities,
there are 43 configurations of singularities on  irreducible tame
torus curves.
The result in \cite{Pho} is valid for non-tame sextics of torus
type
as
  the sub-configurations of the inner singularities on sextics of torus type.
We call them {\em the inner configuration}.
 In this paper, we 
complete the classification of configurations 
of the singularities on irreducible sextics of torus type.

This paper is composed as follows. In \S 2,
we give the list of topological types for 
 outer singularities and explain basic degenerations
among  singularities.
In \S 3, we study  possible outer configurations of singularities.
We start from a given 
inner configuration, and
we determine  the  possible  singularities which can be inserted
outside of the conic $C_2$.
We prove that there exist 121 configurations of  singularities of non-tame sextics
among which there exist 21  maximal configurations (Theorem \ref{classification},
Corollary \ref{Max-config}).
In \S 4,
 we introduce the notions of {\em a distinguished configuration moduli
 space}
 and  
{\em a reduced configuration moduli space} and    {\em a minimal
moduli slice}.
Minimal moduli slices are very convenient for the topological study of 
plane curves.
We   prove that the dimension of a minimal slice is equal to the expected minimal   dimension (Theorem \ref{transversality}).
In \S 5, we give normal forms for the maximal configurations.
In this process, we show that the moduli spaces of certain configurations
are not irreducible
% and  two components
but their minimal slices
 have dimension zero and they have normal forms 
 which  are mutually interchangeable
by a Galois action. However it is not clear if they are isomorphic in the
classical topology. See Proposition \ref{Galois1} and 
Proposition \ref{Galois2}.
We also prove that there exist 5 pairs and  a triple of configurations for  which 
the  moduli spaces  are identical if we ignore
the distinction of  inner and outer singularities (Theorem
\ref{same-moduli}).

For  reduced non-irreducible sextics of torus type,
we will study  their configurations in \cite{Reduced}.
%\cite{Reduced}.

%In \cite{Atlas}, the first author will show that three exists only three
%possibilities for the Alexander polynomial.

%\input NT-2
\section{Inner and outer singularities}
\subsection{Inner and outer singularities}
Let $C$ be an irreducible sextics defined by 
$f(x,y)=0$ where $f=f_2^3+f_3^2$ and  $f_i(x,y)$ is a polynomial of
degree $i$ for $i=2,3$. Here $(x,y)$
is the  affine coordinates $x=X/Z,y=Y/Z$.
% For the brevity's sake, we assume that 
%the line at infinity $Z=0$ is transverse to $C$. 
Let $C_2,C_3$ be the 
conic and the cubic defined by $f_2=0$ and $f_3=0$ respectively.
We assume that the line at infinity is not a component of any of
$C_2, C_3$  and $C$. 
Let $P$ be a singular point of $C$.
A singular point  $P$ of $C$ is called an {\em inner} singularity
(respectively an {\em outer} singularity)
if $P$ is on the intersection $C_2\cap C_3$ (resp. $P\notin C_2\cap C_3$) with respect to 
the torus decomposition (\ref{fixed-torus-decomp}).
We will see later that the notion of {\em inner} or  {\em outer} singularity
 depends on
the choice of a torus expression.
In \cite{Pho}, second author classified inner singularities.
Simple singularities which appear as  singularities on sextics of torus
 type
 are 
% the following:
$ A_2,A_5,A_8,A_{11},A_{14}, A_{17}, E_6, D_4, D_5$.
%Observe that they are all simple singularities.
We use the  following normal forms. % of simple singularities.
% which appear on sextics of torus type:
\[A_n:  y^2+x^{n+1}=0,\quad
E_6: y^3+x^4=0,\quad
D_k: y^2x+ x^{k-1}
\]
%Simple singularities which we encounter as  singularities of  sextics
%of torus curves are 
Non-simple inner singularities on irreducible sextics of torus type  are
the  following (\cite{Oka-Pho1}):
$B_{3,2j}, j=3,4,5 $, $B_{4,6}$, $C_{3,k}, k=7,8,9,12,15$,
$C_{6,6}, C_{6,9}, C_{9,9}$ and $Sp_1$ where
\begin{eqnarray*}
\begin{cases}
B_{p,q}:  \: y^p+x^q=0  \: (\text{Brieskorn-Pham type})\\
C_{p,q}:  \: y^p+x^q+x^2y^2=0,\quad
\frac 2p+\frac 2q\,<\, 1\\
%D_{4,7}:  \: y^4+x^3y^2+x^7=0\\
Sp_1:  \: (y^2-x^3)^2+(xy)^3=0\\
%Sp_2:  \: (y^2-x^3)^2-y^6=0
\end{cases}
\end{eqnarray*}

Note that $B_{3,3}=D_4$.
 For outer singularities,
 a direct computation gives  the following.
\begin{Proposition}\label{outer-singularity}Assume that 
$C$ is an irreducible sextics of torus type and $P\in C$ 
is an outer singularity with  multiplicity $m$. Then
$m\le 3$ and 
the local  topological type $(C,P)$ is a simple singularity and it takes
 one of the following.
\begin{enumerate}
\item If $m=2$, $(C,P)$ is equivalent to one of $ A_1,A_2,\dots,A_5$.
\item If $m=3$, $(C,P)$ is equivalent to  one of $D_4, D_5, E_6$.
\end{enumerate}

\end{Proposition}
\begin{Remark}
The assertion is true for reduced sextics of torus type without 
the assumption of irreducibility.
\end{Remark}
{\em Proof.}
First observe that the sum of Milnor numbers of inner singularities
is bounded   by 12 from below, as the generic  sextics of torus 
type has 6 $A_2$-singularities. By this observation
and by the  genus formula (see \S 3),
the sum of Milnor numbers of outer  singularities is less than or equal
to 
$20-12=8$. By the lower semi-continuity of Milnor number,
the Milnor number of $(C,P)$ is greater than or equal to $(m-1)^2$,
where  $m$  is the multiplicity of $C$ at $P$.
Thus  we get $m\le 3$. The rest of the assertion is proved by an easy
computation.
We may assume that $P=O$, where $O$ is  the origin.
The generic form of $f_2,f_3$ are given as
\begin{eqnarray}\label{normal-form}
 \begin{cases}
& f_2 \left( x,y \right) :={\it a_{02}}\,{y}^{2}+ \left( {\it a_{11}}\,x+{\it a_{01}
} \right) y+{\it a_{20}}\,{x}^{2}+{\it a_{10}}\,x+{\it a_{00}}
\\
&{\it f_3} \left( x,y \right) :={\it b_{03}}\,{y}^{3}+ \left( {\it b_{12}}\,x
+{\it b_{02}} \right) {y}^{2}+ \left( {\it b_{21}}\,{x}^{2}+{\it b_{11}}\,x+{
\it b_{01}} \right) y\\
&+{\it b_{30}}\,{x}^{3}+{\it b_{20}}\,{x}^{2}+{\it b_{10}}\,x+
{\it b_{00}}
\end{cases}
\end{eqnarray}
The condition $P\in C$ and $P\notin C_2$ says that 
$a_{00}=-t^2, b_{00}=-t^3$ for some $t\in \bfC^*$.
Using the condition $f_x (O) =f_y(O)=0$
where $f_x,f_y$  are partial derivatives in $x$ and $y$ respectively,
we eliminate 
coefficients $b_{01}$ and $b_{10}$ as
\begin{eqnarray*}
\mathit{b_{01}} := {  \frac {3}{2}} \,\mathit{t_0}\,
\mathit{a_{01}},\quad
\mathit{b_{10}} := {  \frac {3}{2}} \,\mathit{t_0}\,
\mathit{a_{10}}
\end{eqnarray*}
\noindent
We denote the Newton principal part of $f$ by $NPP(f,x,y)$.
Assume that $m=2$.
Then   $(C,O) = A_1$ generically.
By the action of $\GL(3,\bfC)$, we can assume that the tangent
direction of $(C,O)$ is given by $y=0$.
The degeneration $A_1\to A_2$ is given by  putting
 $f_{xx}(O)=f_{xy}(O)=0$. A direct computation   shows that
the equivalent class  $(C,O)$ can be 
 $A_k$ for $k\le 5$. For example,
to make the degeneration $A_2\to A_3$, we put the coefficient
of $x^3$ in $NPP(f,x,y)$ to be zero.  Then $NPP(f,x,y)$
takes the form $c_2 y^2+c_1 y x^2+ c_0 x^4$ with $c_2\ne 0$ as $m=2$.
The degeneration $A_3\to A_4$ takes place when the discriminant of the
above polynomial vanishes.
Then we take a new coordinate system $(x,y_1)$ so that 
$c_2 y^2+c_1 y x^2+ c_0 x^4=c_2 y_1^2$.
Then we repeat a  similar argument. 
We can see that $A_5\to A_6$ makes $f$ to be divisible by $y^2$ by an easy
computation.
%We can  also confirm this  as $A_6+6A_2$ does not appear  as a lattice
%or sublattice in the list of 
%Yang \cite{Yang}.

\noindent
Assume  that $m=3$.
Generically this gives $(C,O)\cong D_4$.
Assume that the 3-jet is degenerated. 
We may assume 
(by a linear change of coordinates) that 
 the tangent cone is defined by $y^2x$ or $y^3$ corresponding either
the number of the components in the  tangent cone is 2 or 1.
Assume that it is given by $y^2x=0$.
Thus the Newton principal part of $f$ is given by
\begin{multline*}
-{\frac {3}{64}}\,{\frac { \left( {{\it a_{01}}}^{2}+4\,{\it a_{02}}\,{{\it 
t_0}}^{2} \right) ^{2}{y}^{4}}{{{\it t_0}}^{2}}}-1/8\, \left( 16\,{{\it 
t_0}}^{3}{\it b_{12}}+12\,{{\it t_0}}^{2}{\it a_{11}}\,{\it a_{01}}+12\,{{\it t_0}
}^{2}{\it a_{02}}\,{\it a_{10}}+3\,{{\it a_{01}}}^{2}{\it a_{10}} \right) x{y}^{2}
\\
-{\frac {3}{64}}\,{\frac { \left( 4\,{\it a_{20}}\,{{\it t_0}}^{2}+{{\it 
a_{10}}}^{2} \right) ^{2}{x}^{4}}{{{\it t_0}}^{2}}}
\end{multline*}
Thus $(C,O)\cong D_5$. Further we  observe by a direct computation that 
$a_{10}^2+4 t_0^2 a_{20}=0$ makes $f$ reducible. Thus no $D_k$ ($k\ge 6$) 
appears.
If the  tangent cone is given by $y^3=0$, a similar argument shows that 
the only possible singularity $(C,O)$ is $E_6$.
\qed

\subsection{Degenerations on sextics of torus type}
We consider the basic degenerations among singularities.
First, the possibility of the degeneration of  outer singularities 
under fixing the inner singularities is as usual:
% Here is the possibilities:
%\[\text{Outer degenerations}:\quad
 $A_1\to A_2\to A_3\to A_4\to A_5\to E_6$ and 
$ D_4\to D_5\to E_6$.
Of course, some of the above singularities does not exist when
the inner configuration is very restrictive (i.e., far from the generic
one
$=[6A_2]$).

The degenerations of inner singularities are studied in \cite{Pho}:
%The degenerations among simple singularities are 
%\[\text{Inner degenerations}:\quad
 $A_2\to A_5\to A_8\to A_{11}\to A_{14}\to A_{17}$ and 
$A_5\to E_6$.
The degeneration of an outer singularity into an inner singularity is
described by
the following.
\begin{Proposition}
\begin{enumerate}
\item An outer $A_1$ and two inner $A_2$'s degenerate into an $E_6$.
\item An outer $A_2$ and three inner $A_2$'s degenerate into a $B_{3,6}$.
\item  An outer $A_3$ and three inner $A_2$'s degenerate into a
      $C_{3,7}$.
\item  An outer $A_4$ and three inner $A_2$'s degenerate into a
      $C_{3,8}$.
\item  An outer $A_5$ and three inner $A_2$'s degenerate into a
      $C_{3,9}$.
\end{enumerate}
There are no other degenerations.
\end{Proposition}
{\em Proof.}
The proof is computational.
We show the first two degenerations in detail and leave the other cases to the
reader.
We start from the normal form $f=f_2^3+f_3^2$ where
$f_2,f_3$ are given as in (\ref{normal-form}).
We assume that $C$ has a node at O which is not on the conic $C_2$.
Putting $f_2(0,0)=-t_0^2$ and $f_3(0,0)=-t_0^3$ for some $t_0\in
\bfC^*$,
 we get the normal form:
\begin{eqnarray*}
\begin{cases}
&{\it f_2} \left( x,y \right) ={\it a_{02}}\,{y}^{2}+ \left( {\it a_{11}}\,x+
{\it a_{01}} \right) y+{\it a_{20}}\,{x}^{2}+{\it a_{10}}\,x-{{\it t_0}}^{2}
\\
&{\it f_3} \left( x,y \right) ={\it b_{03}}\,{y}^{3}+ \left( {\it b_{12}}\,x+
{\it b_{02}} \right) {y}^{2}+ \left( {\it b_{21}}\,{x}^{2}+{\it b_{11}}\,x+3/2
\,{\it t_0}\,{\it a_{01}} \right) y\\
&+{\it b_{30}}\,{x}^{3}+{\it b_{20}}\,{x}^{2}+
3/2\,{\it t_0}\,{\it a_{10}}\,x-{{\it t_0}}^{3}
\end{cases}
\end{eqnarray*}
We can put $t_0\to 0$ in this form to see  that
two inner $A_2$  singularities are used to create a $E_6$ singularity:
 $2A_2+A_1\to E_6$.
Note that as $f(O)=-t_0^2$, $t_0\to 0$ implies 
the conic $C_2$ approaches to $O$  so that $O$ becomes 
an inner singularity for $t_0=0$.
  To check the degeneration of inner $A_2$  singularities,
we can look at the resultant $R(f_2,f_3,y)$ of $f_2$ and $f_3$ and find that $x=0$ 
has a multiplicity two in $R=0$.

Next we consider that the case $(C,O)=A_2$. We may assume that 
the tangent cone at $O$ is given by $y=0$.
The corresponding normal form is given by
\begin{eqnarray*}
\lefteqn{{\it f_2} \left( x,y \right) ={\it a_{02}}\,{y}^{2}+ \left( {\it a_{11}}\,x+
{\it a_{01}} \right) y+{\it a_{20}}\,{x}^{2}+{\it A_{10}}\,{\it t_0}\,x-{{\it t_0}}^{2}}
\\
\lefteqn{{\it f_3} \left( x,y \right) ={\it b_{03}}\,{y}^{3}+ \left( {\it b_{12}}\,x+
{\it b_{02}} \right) {y}^{2}+ \left( {\it b_{21}}\,{x}^{2}+ \left( -3/4\,{
\it a_{01}}\,{\it A_{10}}+3/2\,{\it a_{11}}\,{\it t_0} \right) x\right )y}\\
&+3/2\,{\it t_0}\,{\it a_{01}} y+{\it b_{30}}\,{x}^{3}-3/8\,{\it t_0}\, \left( {{\it 
A_{10}}}^{2}-4\,{\it a_{20}} \right) {x}^{2}+3/2\,{{\it t_0}}^{2}{\it A_{10}}\,x
-{{\it t_0}}^{3}
\end{eqnarray*}
Here we have substituted  $a_{10}=A_{10} t_0$
%to kill the denominator of the coefficients of  $f$ 
so that we can easily see the limit $\lim_{t_0\to 0}f_i(x,y)$.
We can see easily 
$(C,O)\to B_{3,6}$.
We observe also that the cubic $C_3$  has a node at $O$ as the limit $t_0=0$
and the intersection multiplicity of $C_2$ and $C_3$ at $O$ is 3.
See \cite{Pho} for the degeneration $B_{3,6}\to C_{3,7}\to C_{3,8}\to 
C_{3,9}$. 

For $(C,O)=A_3$,
 the normal form is given as follows and the assertion is easily checked
 by putting $t_0=0$.

%\noindent
%{\bf Case} $(C,O)=A_3$:
\begin{multline*}
{\it f_2} \left( x,y \right) ={\it a_{02}}\,{y}^{2}+ \left( {\it a_{11}}\,x+
{\it a_{01}} \right) y+{\it a_{20}}\,{x}^{2}+{\it A_{10}}\,{\it t_0}\,x-{{\it t_0
}}^{2}
\\
{\it f_3} \left( x,y \right) ={\it b_{03}}\,{y}^{3}+ \left( {\it b_{12}}\,x+
{\it b_{02}} \right) {y}^{2}+ \left( {\it b_{21}}\,{x}^{2}+ \left( -3/4\,{
\it a_{01}}\,{\it A_{10}}+3/2\,{\it a_{11}}\,{\it t_0} \right) x\right)y
\qquad\qquad\\
+3/2\,{\it t_0}\,
{\it a_{01} } y-1/16\,{\it A_{10}}\, \left( {{\it A_{10}}}^{2}+12\,{\it 
a_{20}} \right) {x}^{3}-3/8\,{\it t_0}\, \left( {{\it A_{10}}}^{2}-4\,{\it 
a_{20}} \right) {x}^{2}+3/2\,{{\it t_0}}^{2}{\it A_{10}}\,x-{{\it t_0}}^{3}
\end{multline*}
 The other cases is similar.

\subsection{List of  configurations of inner singularities }

For the classification of non-tame configurations, we start from the 
classification of the configurations of  singularities on tame sextics
of torus 
type \cite{Pho}. 
The list of configurations in \cite{Pho} is  valid as the
sub-configuration defined by
 the inner singularities for a sextics which may have  outer singularities.
%The basic idea was to use the geometry of the conic
%$C_2$ and the cubic $C_3$. 
Let $C_2\cap C_3=\{P_1,\dots, P_k\}$.
The i-vector is by definition the k-tuple of integers
given by the intersection numbers $I(C_2,C_3;P_i),~i=1,\dots,k$.  There
exist
43 possible configurations as follows, assuming $C$ is irreducible. Put
$\bfv:=\iv(C)$

\nin
1.  $\bfv = [1, 1, 1, 1, 1, 1]:\quad$% there is a unique configurations:
%t1. 
$t1=[6A_2]$.

\nin
2.   $\bfv = [1, 1, 1, 1, 2]:\quad$
%, two configurations are possible:
%t2. 
$t2=[4\,{A_{2}}, \,{A_{5}}]$,
%t3. 
$t3=[4\,{A_{2}}, \,{E_{6}}]$.

\noindent
3.   $\bfv=[1,1,2,2]:\quad $
%, there exist 3 configurations:
%t4. 
$t4=[2\,{A_{2}}, \,2\,{A_{5}}]$,
%t5. 
$t5=[2\,{A_{2}}, \,{A_{5}}, \,{E_{6}}]$,
$t6=[2\,{A_{2}}, \,2\,{E_{6}}]$.

\noindent
4.   $\bfv=[1,1,1,3] :\quad $
%, we have 5 configurations:
 $t7=[3\,{A_{2}}, \,{A_{8}}]$,
 $t8=[3\,{A_{2}}, \,{B_{3, \,6}}]$,
$t9=[3\,{A_{2}}, \,{C_{3, \,7}}]$,
 $t10=[3\,{A_{2}}, \,{C_{3, \,8}}]$,
$t11=[3\,{A_{2}}, \,{C_{3, \,9}}]$

\noindent
5.   $\bfv=[2,2,2]:\quad $
%, there exist 4 configurations:
$t12=[3\,{A_{5}}]$,
$t13=[2\,{A_{5}}, \,{E_{6}}]$,
$t14=[A_5,2E_6]$,
$t15=[3\,{E_{6}}]$.

\noindent
6.    $\bfv=[1,2,3]:\quad $
%, there exist 6 configurations:
$t16=[{A_{2}}, \,{A_{5}}, \,{A_{8}}]$,
 $t17=[{A_{2}}, \,{A_{5}}, \,{B_{3, \,6}}]$,
$t18=[{A_{2}}, \,{A_{5}}, \,{C_{3, \,7}}]$,
$t19=[{A_{2}}, \,{E_{6}}, \,{A_{8}}]$,
$t20=[{A_{2}}, \,{E_{6}}, \,{B_{3, \,6}}]$,
$t21=[{A_{2}}, \,{E_{6}}, \,{C_{3, \,7}}]$,

\noindent
7.    $\bfv= [1, 1, 4]:\quad $
%, there exist 5 configurations:
%\indent
$t22=[2\,{A_{2}}, \,{A_{11}}]$,
$t23=[2\,{A_{2}}, \,{C_{3, \,9}^\natural}]$,
$t24=[2\,{A_{2}}, \,{B_{3, \,8}}]$,
$t25=[2 A_2,\, C_{6,6}]$,
$t26=[2A_2, B_{4,6}]$.

\noindent
8.   $\bfv=[3,3]:\quad $
%, we have:
$t27=[2\,{A_{8}}]$,
$t28=[{A_{8}}, \,{B_{3, \,6}}]$,
$t29=[{A_{8}}, \,{C_{3, \,7}}]$,

\noindent
9.    $\bfv=[2,4]:\quad $
%, we have:
%\indent
$t30=[{A_{5}}, \,{A_{11}}]$,
$t31=[{A_{5}}, \,{C_{3, \,9}^\natural}]$,
$t32=[{A_{5}}, \,{B_{3, \,8}}]$,
$t33=[{E_{6}}, \,{A_{11}}]$,
$t34=[{E_{6}}, \,{C_{3, \,9}^\natural}]$,
%\indent
$t35=[{E_{6}}, \,{B_{3, \,8}}]$

\noindent
10.    $\bfv=[1,5]:\quad $
%, we have:
$t36=[{A_{2}}, \,{A_{14}}]$,
$t37=[{A_{2}}, \,{C_{3, \,12}}]$,
$t38=[{A_{2}}, \,{B_{3, \,10}}]$,
$t39=[{A_{2}}, \,{C_{6, \,9}}]$,
$t40=[{A_{2}}, \,{\mathit{Sp}_{1}}]$.

\noindent
11.   $\bfv=[6]:\quad  $
%,  we have  three possibilities:
%\indent
$t41=[{A_{17}}]$,
$t42=[C_{3,15}]$, $t43= [C_{9,9}]$.

\section{Configurations of non-tame sextics}
\subsection{Genus  admissible configurations}
For the classification,  we consider
% we list the possible configurations which satisfies
two inequalities by the positivity of the genus formula:
\begin{eqnarray}\label{genus-formula}
 g(C)=\frac{(d-1)(d-2)}2 -\sum_{P\in \Si(C)}\de(C,P)\ge 0
\end{eqnarray}
and by the positivity of the class number $n^*(C)$:
\begin{eqnarray}\label{class-number}
 n^*(C)=d(d-1)-\sum_{P\in \Si(C)} (\mu(C,P)+m(C,P)-1)\,\ge\, 0
\end{eqnarray}
Here $d=\degree(C)$, % ($d=6$ in our situation),
  $\Si(C)$ is the set of singular points of $C$ and 
$\de(C,P)$ is the $\de$-genus of $C$ at $P$  which is equal to
$\frac 12 (\mu(C,P)+r(C,P)-1)$ with $r(C,P)$ being the number of local
irreducible components at $P$ (see Milnor \cite{Milnor}).
The class number $n^*(C)$ of $C$ is defined by the degree of 
the dual curve $C^*$  where  $m(C,P)$ is the multiplicity of $C$ at 
$P$. 
See \cite{NambaBook,Oka-sextics} for the class number formula
(\ref{class-number}).

A configuration $\Si$ is called a {\it genus-admissible}
%(respectively {\it existing}) 
if   the  genus  and   the class 
number given by the above formulae (\ref{genus-formula}),
(\ref{class-number})  are non-negative.
%(resp. if there exists an irreducible sextics of torus type
%with $\Si$ as the configuration of the singularities).

There exist 145 configurations which satisfy those inequalities. See
Table
1
$\sim$ 5 in the end of this paper.
  In the  list,
% of  145 genus-admissible configurations
%of non-tame torus sextics, 
the first  bracket
shows the  configuration of the  inner singularities and the second is 
that of  the outer
 singularities.
 For example, $[[6A_2], [3A_1]]$  shows that $C$ has 6 $A_2$ as inner
 singularities
 and 3 $A_1$ as outer singularities.
% The invariant vector 
 %is defined by 
The vector $(g(C), \mu^*(C),n^*(C),i(C))$ denotes the invariants of $C$,
where
 $g(C)$ is the genus of the normalization, $\mu^*(C)$ is the sum of
 Milnor numbers at singular points,
  $n^*(C)$ is the class number and $i(C)$ is defined by
$3d(d-2)-\sum_{P}\de(P)$ which is the
 number of flex points on $C$.
For the calculation of $\de(P)$, we refer
Oka \cite{Oka-sextics}. 
(In Corollary 12 of \cite{Oka-sextics}, 
there is a trivial mistake. The correct formula
is $\bar\de(A_{2p-1})=6p$ for any $p$ which follows from Theorem
10, \cite{Oka-sextics}.)

\subsection{Existing configurations}
The main problem  is how to know those  configurations which do exist
and
which do
 not exist in the list of Table 1 $\sim$ 5 in subsection \ref{Table}. 
\begin{Theorem} \label{classification}
The possible configurations  of singularities of
irreducible sextics of torus type with at least one outer 
singularity
are given by  Table1 $\sim$ Table 5 in the last subsection \ref{Table}.
There are 24 configurations  in the table  which do not exist
(they are marked `No') and the other 121
 configurations  exist. 
% Among them, there exists 21 
%maximal configurations.

Combining the list of the configuration of tame  sextics of torus type,
there exist 164 configurations
on irreducible sextics of torus type.
\end{Theorem}

 The  column of the table 
``Existence ?'' provides the informations about existence and
 non-existence and typical degenerations. ``No'' implies the
 corresponding configuration does not exist. ``Max'' implies that
 the configuration is maximal among irreducible 
sextics of torus type. The arrow shows a possible
 degeneration.
The last column gives the expected minimal moduli slice dimension,
which is defined in \S 3.
\begin{Corollary}
The fundamental group $\pi_1(\bfP^2-C)$ of the complement of a sextics
$C$ with
a configuration corresponding to one of
 the following 
%numbers in the  list of non-tame  configurations
 is isomorphic to 
$\bfZ_2*\bfZ_3$ by \cite{Oka-Pho1}.
\[nt\,j,\,j=
 1,2,3,4, 5,
19,25,26,27, 33,43,44,45,54,61,68,72,73,74,90, 92
\]
\end{Corollary}
%\begin{Remark}
%There are 5 pairs  and a triple of configurations for which the
%configurations are identical if we ignore the distinction of 
%the inner and outer singularities.
%
%They gives actually the same moduli spaces.
%
\begin{Corollary}\label{Max-config}
There exist 21 maximal configurations on non-tame
sextics of torus type:
\begin{eqnarray*}
%\lefteqn{
&& nt23= [[6\,A_2],[3A_2]],~nt32=[[4\,A_2,A_5],[E_6]],~
nt47=[[4\,A_2,E_6],[A_5]]\\
&&
  nt64=[[2\,A_2,A_5,E_6],[A_4]],
 nt67=[[2\,A_2,A_5,E_6],[2\,A_2]], nt70=[[2\, A_2,2\,E_6],[A_3]]\\
&&
 nt78=[[3\,A_2,A_8],[D_5]],
  nt83=[[3\,A_2,A_8],[A_1,A_4]],
nt91=[[3\,A_2,B_{3,6}],[A_2]]\\
&&
 nt99=[[A_5,2E_6],[A_2]],
 nt100=[[3E_6],[A_1]],
nt104=[[A_2,A_5,A_8],[A_4]]
\\
 &&
nt110=[[A_2,E_6,A_8],[A_3]],
nt113=[[A_2,E_6,A_8],[A_1,A_2]], nt118=[[2\,A_2,A_{11}],[A_4]]
\\
&& nt123=[[2\,A_2,C_{3,9}^\natural],[A_2]], nt128=[[2A_8],[A_3]], 
nt136=[[E_6,A_{11}],[A_2]]\\
&& 
nt139=[[A_2,A_{14}],[A_3]],
nt142=[[A_2,A_{14}],[A_1,A_2]], nt145=[[A_{17}],[A_2]]
\end{eqnarray*}
\end{Corollary}

In the table, $C_{3,9}$ and $C_{3,9}^{\natural}$ are
topologically isomorphic but 
they are distinguished by  $\iota=3$ and $4$ respectively. See Pho \cite{Pho}. 

%\end{Remark}
%\input Non-existence

\vspace{.5cm}
\subsection { Proof of the non-existence of 24 configurations}
In this subsection, we prove  the non-existence of
the configurations
nt j, $j=$14,16, 17,
18, 38, 39, 48, 76, 79, 84, 85, 86, 89, 93,
 107, 111, 114, 121, 125, 129, 131, 132, 140, 143 in Table $1\sim 5$.
It is well-known that the total sum of the Milnor numbers $\mu^*$
on sextics is bounded by 19 if the singularities are all simple
(\cite{Horikawa:1975},\cite{Shioda-Inose}).
Thus the configurations nt79, nt86, nt111, nt114, nt129, nt132, nt140, nt143 do not exist.

Another powerful tool is to consider the dual curves.
We know that the dual singularities of $A_k,k\ge 3$,
$C_{3,p},p\ge 7$ and  $B_{3,q},q\ge 6$ are generically
isomorphic to  themselves
 \cite{Oka-sextics}.  If the singularity is not generic,
the dual singularity has a  bigger Milnor number. 
The  singularity  $B_{3,3}$ corresponds
to a tri-tangent
 line in  the dual curve $C^*$.
By Bezout theorem,  a tri-tangent line does not exist for curves of degree
$\le 5$.
Thus the existence of  $B_{3,3}$ implies $n^*(C)\ge 6$.

The non-existence of the configurations nt14, nt16, nt17, nt18, nt38,
nt39,
nt48, nt89, nt93,
and nt125
 can be proved by taking the dual curve information into
consideration.
For example, consider the configuration  nt14 $=[[6A_2],[A_1,B_{3,3}]]$.
If such a curve $C$ exists, the dual curve $C^*$ has degree 4, which
 is impossible.
% as a quartic does not have a tri-tangent line.
Next we show that the configurations nt16 $\sim$ nt18 do not exist.
Assume a curve $C$ with 
the configuration nt16$=[[6A_2],[A_2,A_3]]$ for example.
 Then  the dual curve $C^*$ has
degree
5 and $C^*$ has $A_3, 4A_2$ as singularities. By the class
formula,  the dual curve $C^{**}=(C^*)^*$ have degree 4
which is absurd. The other two can be eliminated in the same discussion.

 For nt93, we use the fact that the dual singularity of
$C_{3,7}$
is again $C_{3,7}$ ( \cite{Oka-Pho1}). Assume that there exists a
sextics $C$
with configuration nt93. Then the dual curve $C^*$ is a quintic with 
$C_{3,7}$ and $A_2$. Then by the Pl\"ucker formula,
this is ridiculous as $\de(C_{3,7})=6$.
Suppose that a sextics with  the configuration nt125 exists.
Then the dual curve have degree 5 and $B_{3,8}$ as a singularity.
However the total sum of the Milnor numbers on an irreducible quintic
is bounded by 12, a contradiction.
 The other configurations are treated in a similar way.

The configurations nt76, nt85, nt107, nt121 and  nt131 do not exist as they are not
in the list   of Yang table
 \cite{Yang}. The non-existence of these configurations can be also checked by a direct
 maple computation. The non-existence of nt84 has to be checked by a 
direct
 computation.
\qed
\begin{Remark}
We remark here that a configuration in the list of Yang does not
 necessarily  exist as a configuration of a sextics of torus type.
There are also a certain configurations with only simple singularities
 which
is not a sublattice of a lattice of maximal rank  in Yang's list.
\end{Remark}

\section{Moduli spaces}

%Let $C$ be a reduced plane curve and let $P_1,\dots, P_k$ be the
%singular points.
%$\{(C,P_i); i=1,\dots,k\}$ is referred as the configuration of the singularities.
%
\subsection{Distinguished configuration moduli and reduced configuration moduli}
%We introduce several moduli spaces which we are interested in.
Let $\Si_1,\Si_2$ be  configurations of singularities. 
In this paper, a configuration is a finite set of 
topological equivalent classes of germs of isolated curve singularities.
We say that 
$\Si:=[\Si_1,\Si_2]$  be a {\em distinguished configuration } on a
sextics of torus type if 
 $\Si_1$ is the configuration of inner singularities 
and $\Si_2$ is the configuration of outer singularities.
%In Table  1 $\sim$ 5 of configuration in Theorem \ref{classification},
%we have used the same notation.
We put
$\Si_{red}:=\Si_1\cup \Si_2$
%This is  the configuration ignoring the distinction
%of inner and outer singularities. 
 and we call $\Si_{red}$ 
a {\em reduced configuration}. We now introduce several moduli spaces
which we consider
in this paper.
First, recall that  spaces of conics  and cubics are 6 and 10 dimensional
respectively. Let $\tilde{\cT}$ be the vector space of dimension 16 which is
% weighted projective space of dimension 15
defined by 
\[
  \tilde{\cT}:=\{\bff=(f_2,f_3);\degree\, f_2=2,~\degree\, f_3=3\}
\]
There is a canonical $\GL(3,\bfC)$-action on $\tilde{\cT}$. The center of
$\GL(3,\bfC)$ is identified with $\bfC^*$. It defines a canonical
weighted homogeneous action on $\tilde{\cT}$ and we introduce an
equivalence relation
$\sim$
by
$(f_2,f_3)\sim (f_2',f_3')\iff f_2'=f_2\, t^2,~f_3'=f_3\, t^3$  for some $t\in
\bfC^*$.  In particular,
$(f_2,f_3)\sim (f_2\, \omega^j,\pm \, f_3)$ for $j=1,2$ where 
$\omega=(\sqrt{3}I-1)/2)$. ( We use the  notation
$I=\sqrt{-1}$.)
Let $\cT$ be the weighted projective space by
the $\bfC^*$-action and let $\pi:\tilde \cT\to \cT$ be the quotient map. Then
$\PGL(3,\bfC)=\GL(3,\bfC)/\bfC^*$ acts on $\cT$.
Each equivalence class $(\bff)$ defines a sextics of torus type
$C(\bff)$ defined by $f_2(x,y)^3+f_3(x,y)^2=0$. 
We put
\[
\Si(\bff)_{in}:=\{(C(\bff),P_i); f_2(P_i)=0 \},\quad
\Si(\bff)_{out}:=\{(C(\bff),P_i);f_2(P_i)\ne 0\}
\]
where  $\{P_1,\dots, P_k\}$
are  the singular points of $C(\bff)$ and 
 $(C(\bff),P_i)$ is the topological equivalent class of
the germ   at $P_i$.
%We put also $\Si(\bff)=[\Si(\bff)_{in},\Si(\bff)_{out}]$.
Let $\Si=[\Si_1,\Si_2]$ be a distinguished configuration.
The {\em distinguished configuration moduli} $\cM(\Si)\subset \cT$ is defined by the quotient
$\tilde{\cM}(\Si)/\bfC^*$
\[
\tilde{\cM}(\Si):=\{\bff=(f_2,f_3)\in \tilde{\cT}; \Si(\bff)_{in}=\Si_1,~\Si(\bff)_{out}=\Si_2 \}
\]
The space of sextics, denoted by $\tilde{\cS}$, is a vector space of
dimension 28
and its quotient by the homogeneous $\bfC^*$-action is denoted by $\cS$.
 There exist a canonical $\GL(3,\bfC)$-equivariant 
mapping $\tilde{\psi}_{red}: \tilde{\cT}\to \tilde{\cS}$ which is defined by
$\tilde{\psi}_{red}(f_2,f_3)=f_2^3+f_3^2$ and it induces a canonical
$\PGL(3,\bfC)$-equivariant
mapping $\psi_{red}:\cT\to \cS$.
Let $\Si_0$ be a reduced configuration.
The {\em reduced configuration moduli} $\cM_{red}(\Si_{0})$
is defined by $\tilde{\cM}_{red}(\Si_{0})/\bfC^*$ where $\bfC^*$-action is the
the scalar
multiplication and
\[
 \tilde{\cM}_{red}(\Si_{0})=\{f\in \tilde{\cS}; \exists
 \Si=[\Si_1,\Si_2],~\Si_0=\Si_1\cup \Si_2,
~\exists\bff\in \cM(\Si), ~\tilde{\psi}_{red}(\bff)=f\}
\]

The map $\psi_{red}: \cM(\Si)\to \cM_{red}(\Si_{red})$ 
is not necessarily injective (see
 Observation \ref{not-injective}). 
% The 
%{\em quotient moduli spaces} 
% $\cM(\Si)/\PGL(3,\bfC)$ or 
%$\cM_{red}(\Si_{red})/\PGL(3,\bfC)$ are not used explicitly in this paper.
%
\begin{Remark}
Let $\bff=(f_2,f_3)\in \tilde{\cM}([\Si_1,\Si_2])$ and assume that
 $f_2^3+f_3^2=0$
is an irreducible sextics and assume that $\Si_2$ is not empty.
Consider the family of sextics
$C_t: t\,f_2^3(x,y)+f_3(x,y)^2=0$. By Bertini theorem, 
for a generic $t\ne 0$, $C_t$ has only inner singularities and 
 $\Si(C_t)=\Si_1'$, where simple singularities in $\Si_1$ are unchanged
in $\Si_1'$
and non-simple  singularities
are replaced by the first generic singularities fixing the singularities
of the conic $f_2=0$ and the cubic $f_3=0$ and their local intersection
 numbers
in Table $A'$ of \cite{Pho}.
For example, inner singularities with a nodal cubic  and a smooth conic, with the
 intersection number $3$, any singularity in the series
$B_{3,6}\to C_{3,7}\to C_{3,8}\to C_{3,9}$ is replaced by $B_{3,6}$.
This is the reason why we need the information of defining polynomials
 $f_2,f_3$,
not only the geometry of $C_2$ and $C_3$.
\end{Remark}

\subsection{Moduli slice and irreducibility} 

A subspace $A\subset \tilde{\cM}(\Si)$ is called  a {\em   moduli 
slice} of $\cM(\Si)$
if its  $\GL(3,\bfC)$-orbit  covers the  whole moduli space   
$\tilde{\cM}(\Si)$ and  $A$ is an algebraic variety. 
%Usually a moduli slice is given by fixing the locations of some of 
 %singularities at explicit points and normalizing 
%a coefficient of $f_2(x,y)$. 
A  moduli slice is called  
{\em   minimal } if the dimension  is minimum. 
As we are  mainly interested in  the topology of 
 the pair $(\bfP^2,C)$ where $C$ is a sextics defined by 
$f_2^3+f_3^2=0$, % for  $(f_2,f_3)\in \cM(\Si)$, 
the important point is the  
connectedness of the moduli. 
Thus we are  interested, not in  the algebraic structure 
of the moduli spaces
but  in the explicit form of a minimal moduli slice,
which we call {\it a normal form}.
Note that the moduli space $\cM(\Si)$  might be irreducible even
 if a minimal  slice $\cA$ is not irreducible.
%For the determination of normal forms,
% without unnecessary parameters, 
%we start from the normal form (\ref{normal-form}) in 2.1.  
%Though it has 15-dimensional, it is better 
%to decrease the parameters without losing generality. 
%

Points $P_1,P_2,\dots, P_k$ in $\bfP^2$  are called  {\em generic} if 
any three of them  are not on a line. Let $P_1,P_2,P_3$ are generic
points and
let $L_i$ be lines through $P_i$, for $i=1,2$. We say $L_i$ is a {\em generic}
line through $P_i$ with respect to $\{P_1,P_2,P_3\} $ if
$L_i$ does not
pass through any of other  two points $\{P_j;j\ne i\}$.  Observe that
two set of  generic four points, or of
 generic three points and two generic lines through two of them 
 are transformed each other by
$\PGL(3,\bfC)$-action.
% and also
% any generic three points and two generic lines through two of them can
% be transformed each other  by $\PGL(3,\bfC)$-action.
%We leave the proof of this assertion to the reader.
Note  that  
%it takes two dimension from $\PGL(3,\bfC))$ to put a singularity at a fixed 
%location. That is,
 the dimension of 
the isotropy 
group of a point (respectively a point and a line through it) is
codimension 2
(resp. 3).
As $\dim \PGL(3,\bfC)=8$, we can fix, using the above principle 
either  
 
(a)	location of four singularities at generic positions or  
 
(b)	three  singularities  at generic positions
and two generic tangent cones. 
 
\noindent
This technique is quite useful to compute various normal forms. 
\subsection{Virtual dimension and transversality}  
In general, the dimension of the moduli space of a given  
configuration of  singularities is difficult to   
be computed. However in the space of sextics of torus type, the situation is quite simple.  
Suppose that we are given a sextics defined by (\ref{normal-form}).  
%It has 15 dimensional, as we consider it in the projective space. 
Take a point $ P=(\al,\be)\in \bfC^2$ and consider the condition for $P$
to be a singular point of $C$. For simplicity we assume that  
$P=(0,0)$.  
  
(I) First assume that $P$ to be an inner singularity.  
Let $\si$ be a class of $(C,P)$.   
We define the integer $\icodim(\si)$ by (the number of  
independent conditions  
on the coefficients) $- 2$. 
Here $2$ is the freedom to choose $P$. 
For example, the condition for $P$ to be an inner $A_2$ singularity is  
simply $f_2(P)=f_3(P)=0$.  
%However if $P$ is arbitrary,  
% $P$ can move in two dimensional space. 
% So this does not give any restriction. 
 So $\icodim (A_2)=0$.  
 Assume that $(C,P)\cong A_5$.  
Then the corresponding condition is $f_2(P)=f_3(P)=0$ and the intersection multiplicity   
of $C_2$ and $C_3$  at $P$ is 2. This condition is equivalent to
$(f_{2x}f_{3y}-f_{2y} f_{3x})(P)=0$.
Thus  $\icodim(A_5)=1$. Similarly the condition $(C,P)\cong E_6$ is given by  
$f_2(P)=f_3(P)=0$ and the partial derivatives $f_{3x}$  
and $f_{3y}$  
vanishes at $P$. See Pho \cite{Pho} for the characterization of inner
singularities. Thus we  have $\icodim(E_6)=2$.  
%$\icodim(A_8)=2$, i-$\codim(A_{11}=3$, i-$\codim (A_{14})=4$  
%and i-$\codim(A_{17})=5$.  
  
Let $\iota=I(C_2,C_3;P)$ be the intersection number 
of $C_2$ and $C_3$ at $P$. Similar discussion proves that   
\begin{Proposition}  
For the inner singularities on  sextics of torus type, 
 $\icodim$ is given as follows.

\begin{tabular}{c|c|c|c|c}
 $\icodim$& 1 & 2& 3& 4\\
singularity&$A_5$& $E_6$, $A_8$& 
 $~A_{11},  B_{3,6}$&
\vtop{\hbox{ $A_{14},C_{3,9}^{\natural}$}
\hbox{ $C_{3,7}, C_{6,6}$}}\\
\hline
$\icodim$&  5&6&7&8\\
singularity&\vtop{\hbox{$A_{17}, C_{3,12}$}
\hbox{$B_{3,8},C_{3,8}$}
\hbox{$C_{6,9}, B_{4,6}$}}&
\vtop{\hbox{$ C_{3,15},B_{3,10}$}\hbox{$ C_{3,9}, Sp_1$}
\hbox{$C_{9,9},C_{6,12}$}\hbox{ $D_{4,7}$}}& 
 $B_{3,12}, Sp_2$& 
 $B_{6,6}$\end{tabular}
\end{Proposition}  
The proof is immediate from the above consideration and the   
existence of the degeneration series where each step is codimension one (\cite{Pho}).
The vertical degenerations keep the intersection number 
$\iota$ and it is observed to have codimension one for each arrow  
in \cite{Pho}. The first and the second horizontal sequence are
 induced by increasing $\iota$ by one for each arrow.  
Thus each arrow has codimension one.  
Recall that $P$ is $C_{6,6}$ singularity if both of $C_2$ and $C_3$ has a node at $P$.   
Thus we can easily see that $\icodim(C_{6,6})=4$.   
The degenerations $C_{6,6} \to C_{6,9}\to C_{9,9}$ or $C_{6,9}\to C_{6,12}$ has also codimension one for each arrow   
(\cite {Pho}). 
So the rest of assertion follows immediately from the   
above consideration.  

\[  
\begin{matrix}  
A_5&\to& A_8 &\to&  A_{11} &           \to & A_{14} &\to &  A_{17}\\  
\mapdown{}&&   \mapdown{}&&\mapdown{}&&\mapdown{}&&\mapdown{}\\  
        E_6&\to& B_{3,6}&\to&C_{3,9}^{\natural}&\to & C_{3,12}&\to &C_{3,15}\\  
&&\mapdown{}&&\mapdown{}&&\mapdown{}&&\mapdown{}\\  
&&              C_{3,7}&\to&B_{3,8}&\to&B_{3,10}&\to& B_{3,12}\\  
&&\mapdown{}&&&&&&\\  
&&C_{3,8}&&C_{6,6}&\to&C_{6,9}&\text{\vtop{\hbox{$\to$}\hbox{$\to$}}}& \text{\vtop{\hbox{$C_{9,9}$}\hbox{$ C_{6,12}$}}}\\  
&&\mapdown{}&&\mapdown{}&&\mapdown{}&&\\  
&&C_{3,9}&&B_{4,6}&\to&Sp_1&\to&Sp_2\\  
 &&&&            \mapdown{}&&&&\\  
&&&&D_{4,7}&&&&  
\end{matrix}\]   
  
(II) Now we assume that $P$ is an outer singularity. This means $f_2(P)\ne 0$. Let $\si$ be the topological equivalent class   
of $(C,P)$.  
We define the integer $\ocodim(P)$ by the number of conditions 
 on the space of coefficients of $f$ minus 2.  
By  the argument in the proof of Proposition \ref{outer-singularity},  we can easily see  
that  
\begin{Proposition} For an outer singularity on sextics of torus type, we have  
$\ocodim(A_i)=i,~i=1,\dots,5$ and $\ocodim(D_i)=i,~i=4,5$ and 
 $\ocodim(E_6)=6$. 
Thus in all cases, $\ocodim(\si)$ is equal to the Milnor number. 
\end{Proposition}  
{\em Proof.} For $A_1$, we need three condition $f(P)= f_x(P)= 
f_y(P)=0$.  
Here $f_x,f_y$ are partial derivatives. 
Thus $\ocodim(A_1)=3-2=1$. The other assertion follows from the   
basic degeneration series of codimension one:  
\[A_1\to A_2\to A_3\to A_4 \to A_5,\quad B_{3,3}=D_4\to D_5\to E_6  
\]  
Note that  $B_{3,3}$ singularity is defined by 6 equations, $f(P)=f_x(P)=f_y(P)=f_{x,x}(P)=f_{xy}(P)=f_{yy}(P)=0$.  
Thus we have $\ocodim (B_{3,3})=4$ and other assertion follows from the above degeneration sequence.  
  
For a given configuration $\Si=[\Si_1,\Si_2]$ on sextics of torus
 type,
% where the singularities in $\Si$ has a mark of inner or outer singularity.  
we define the {\em expected minimal moduli slice dimension}, denoted by   
$\emdim (\Si)$ by the integer  
$$\emdim (\Si):=16-\sum_{\si\in \Si_1}\icodim(\si)-\sum_{\si\in \Si_2}\ocodim(\si)-9$$  
Here 16 is the dimension of sextics of torus type and $9$ is the dimension of $\GL(3,\bfC)$.  
On the other hand, we denote the dimension of minimal moduli slice of the decomposition
moduli $\cM(\Si)$ by $\msdim(\Si)$.   
By the above definition, it is obvious that   
%the minimal moduli slice dimension 
\[
 \msdim(\Si)\ge\emdim(\Si)\]
if it is 
not empty.  
We say that $\Si$ has a {\em transverse moduli slice} or the moduli 
 $\cM(\Si)$ is {\em transversal}
 if $\msdim(\Si)=\emdim(\Si)$. 
%\end{Definition} 
%
%By a direct computation, we can show that   
\begin{Theorem}\label{transversality} 
For  any configuration  $\Si$ of sextics of torus type, there exists
a component $\cM_0$ of $\cM(\Si)$ whic is transvesal.
%\newline\indent 
% $\cM_0$ has the expected minimal  moduli slice dimension.  
\end{Theorem}  
It is probably true that $\cM(\Si)$ is transverse for any $\Si$ but 
we do not want to check this assertion for 121 cases.
The proof of the above weaker assertion
 is reduced to the assertion on maximal configurations
(see the next section) and to the following proposition.  
\begin{Lemma}\label{Red-Maximal} Assume that  $\Si$ degenerates into a maximal 
 configuration  $\Si'$ which has a transverse moduli slice.
% and is irreducible.  
Then the moduli $\cM(\Si)$ has also a component which
is transversal.
% has the expected minimal  moduli slice  dimension. 
\end{Lemma}   
{\em Proof.}   
By the definition, a minimal moduli slice  for $\cM(\Si')$ can be obtained by adding   
$\nu$ equations on the space of coefficients where 
$\nu:=\emdim(\Si)-\emdim(\Si')$.  
Thus we have  
\begin{multline*} 
\msdim(\Si')\ge \msdim(\Si)-(\emdim(\Si)-\emdim(\Si'))\\ 
\ge \msdim(\Si')+(\msdim(\Si)-\emdim(\Si)) \ge \msdim(\Si')
\end{multline*} 
which implies the assertion. \qed

\section{Minimal moduli slices for maximal configurations} 
In this section, we give  normal forms of  minimal moduli slices 
for 
the maximal configurations.
% and some other configurations. 
Using the degeneration argument and Lemma \ref{Red-Maximal}, this guarantees the existence of any
other non-maximal configurations in Table 1 $\sim$ 5 in the subsection \ref{Table}.
We also show that they have transverse minimal 
 moduli slices. 

\vspace{.2cm}
\noindent 
{\bf nt23.}  We consider the minimal  moduli slice of 
 $\cM(\Si_{23})$ with $\Si_{23}=[[6A_2],[3A_2]]$ 
 by the following  minimal slice 
condition: 
 
($\star$) three outer $A_2$'s are  at $P_0:=(0,0)$ and $P_{1}:=(1, 1)$
and $P_3:=(1,-1)$.
The (reduced) tangent cones of $C$ at $(1,\pm 1)$ are given by  $y=\pm 1$ respectively.
 
The calculation is easy. We start from the normal form
$f=f_2^3+f_3^2$ where $f_2,f_3$ are given as in (\ref{normal-form}).
Necessary conditions are
\[
 f_2(P_i)=-t_i^2,\quad f_3(P_i)=-t_i^3,\quad  
 f_x(P_i)= f_y(P_i)=0,~i=0,1,2.
\]
The assumption on the tangential cones gives
$f_{xy}(P_i)= f_{xx}(P_i)=0, \quad i=1,2$.
Solving these equations, we get  the following normal form
with  one free parameter $t:=t_0$.  As $\emdim(\Si_{23})=1$, 
it has a transverse minimal moduli slice.
% whose  normal form is given by $f=f_2^3+f_3^2$ where 

\begin{multline*}
\mathrm{f_2}(x, \,y)=y^{2} - 9\,x^{2}\,t^{2} - 3\,x^{2} + 6\,t^{2
}\,x + 2\,x - t^{2}
\\
\mathrm{f_3}(x, \,y)={\displaystyle \frac {1}{2}} ( - 9\,y^{2}\,t
^{2}\,x - 3\,y^{2}\,x + 3\,y^{2}\,t^{2} + 3\,y^{2} + 3\,x^{3} + 
27\,x^{3}\,t^{2} + 54\,x^{3}\,t^{4} - 3\,x^{2} - 27\,x^{2}\,t^{2}
 \\
\mbox{} - 54\,t^{4}\,x^{2} + 18\,x\,t^{4} + 6\,t^{2}\,x - 2\,t^{4
})/t 
\end{multline*}
As is well-known, the corresponding sextics are the dual of smooth cubics.
 
\vspace{.2cm}
\noindent 
{\bf nt32.} We consider the moduli space $\cM(\Si_{32})$ with 
$\Si_{32}=[[4A_2,A_5],[E_6]]$. 
The irreducibility is easily observed using  the slice condition:
\newline
$(\star)$ an inner 
$A_5$ is  at $(0,0)$ and an outer $E_6$ is at $(0,1)$
with respective
 tangent cones defined by
$y=0$ and $y=1$.

We usually use the $\bfC^*$-action to normalize the coefficient of $y^2$
in $f_2$  to be 1.
The normal forms are given by

%\label{normal32}
\begin{eqnarray*}
{\lefteqn{\mathrm{f_2}(x, \,y)=y^{2} + ( - 1 - \mathit{t_1}^{2})\,y + 
\mathit{a_{02}}\,x^{2}}}\\
\lefteqn{\mathrm{f_3}(x, \,y)={\displaystyle \frac {1}{8 \mathit{t_1}}} \,
(\mathit{t_1}^{4} + 6\,\mathit{t_1}^{2} - 3)
\,y^{3}  + {\displaystyle \frac 
 {1}{8\mathit{t_1}}} \,(6 - 6\,\mathit{t_1}^{4})\,y^{2}}
 \\
&\mbox{} + {\displaystyle \frac  {1}{8\mathit{t_1}}} \,(
 - 6\,\mathit{a_{02}}\,x^{2} - 6\,\mathit{t_1}^{2} + 6\,\mathit{t_1}^{
2}\,\mathit{a_{02}}\,x^{2} - 3\,\mathit{t_1}^{4} - 3)\,y
  + {\displaystyle \frac  {1}{8 \mathit{t_1}}} \,(6\,
\mathit{t_1}^{2}\,\mathit{a_{02}}\,x^{2} + 6\,\mathit{a_{02}}\,x^{2})
\end{eqnarray*}

Observe that $\cM(\Si_{32})$ is irreducible by this expression.
We have used 6 dimension of $\PGL(3,\bfC)$ for the above slice.
To get a minimal slice, we have two more dimension to use, so we can fix
a location of an inner $A_2$. 
Here, we have two choice: either (a) to choose a location  which is
on $\bfQ^2$
or (b) to choose a simple normal form.
The case (a) give as a little complicated normal form. So we choose (b).
We choose $t_1=a_{02}=1$. This can be done by taking 
an inner $A_2$-singularity  at 
$(\al,\be)$ where 
\[
 \al=-\frac 12 \sqrt{6-2I\sqrt{3}},\quad \be=(3+I\sqrt{3})/2.
\]
Note that $\al$ is not well-defined but $\al^2$ is well-defined. This is
enough as $f_2(x,y),f_3(x,y)$ are even in $x$ in the above normal form
%(\ref{normal32})
and the condition implies also $(-\alpha,\beta)$ is another inner $A_2$. 
The corresponding minimal slice has dimension 0, and consists of two
points
and as  the moduli is irreducible, we can take the normal form
\begin{eqnarray}
&f_2(x,y)\, =\, y^2-2\,y+ x^2,\qquad
f_3(x,y)\,=\,(y^3-3\, y+3\, x^2)/2\notag\\
&f(x,y)=( y^2-2\,y+ x^2)^3
+ (y^3-3\, y+3\, x^2)^2/4\label{reduced-32}
\end{eqnarray}
Let $f_{(32)}(x,y)$ be the corresponding sextics.

%input nt47-normal.tex 
\vspace{.2cm}\noindent 
{\bf nt47.} The moduli space $\cM(\Si_{47})$
is irreducible and $\emdim(\Si_{47})=0$ where
 $\Si_{47}:=[[4 A_2, E_6], [A_5]]$.
 This can be checked easily using the slice as in nt32:

$(\star)$  an outer $A_5$ is at $(0,0)$ and an inner $E_6$ is at $(0,1)$
with respective 
 tangent cones given by
$y=0$ and $y=1$.

The corresponding normal form is given as
\begin{eqnarray*}
\lefteqn{\qquad\mathrm{f_2}(x, \,y)=y^{2} + ( - 1 + \mathit{t_0}^{2})\,y - 
\mathit{t_0}^{2} + \mathit{a_{02}}\,x^{2}}\\
&\mathrm{f_3}(x, \,y)=( - {\displaystyle \frac {3}{2}} \,\mathit{
t_0} - {\displaystyle \frac {1}{2}} \,\mathit{t_0}^{3})\,y^{3} + 3
\,y^{2}\,\mathit{t_0} + y\, \left(  \! {\displaystyle \frac {3}{2}
} \,\mathit{t_0}\,( - 1 + \mathit{t_0}^{2}) + {\displaystyle 
\frac  {3}{2\,\mathit{t_0}}} \,\mathit{a_{02}}\,x^{2}  \!  \right) 
 - \mathit{t_0}^{3}
+ {\displaystyle 
\frac {3}{2}} \,\mathit{t_0}\,\mathit{a_{02}}\,x^{2}
\end{eqnarray*}

Thus the irreducibility follows from this expression. Observe also that 
$f_2,\, f_3$ are even in $x$.
 Now we compute the minimal moduli slice with an additional condition,
 an inner $A_2$ at $(\al,\be)$ where $\al,\be$ are  as in nt32.
As a minimal slice, we can take
\begin{eqnarray*}
&\mathrm{f_2}(x, \,y)=y^{2} - {\displaystyle \frac {5}{2}} \,y + 
{\displaystyle \frac {3}{2}}  + {\displaystyle \frac {1}{2}} \,x
^{2}\\
&
\mathrm{f_3}(x, \,y)=\sqrt{6}\, I\,\left ( - {\displaystyle \frac {3}{8}} 
\,y^{3} + {\displaystyle \frac {3}{2}} \,y^{2} + 
{\displaystyle \frac {1}{6}} \,( - {\displaystyle \frac {45}{4
}}  - {\displaystyle \frac {3}{2}} \,x^{2})\,y + {\displaystyle 
\frac {1}{6}} \,({\displaystyle \frac {9}{2}}  + 
{\displaystyle \frac {9}{4}} \,x^{2})\right )
\end{eqnarray*}

Let $f_{(47)}(x,y)$ be the corresponding sextics.
The above normal form  proves $\msdim(\Si_{47})=\emdim(\Si_{47})=0$. 
%$Let $f_{32}(x,y)$  be the sextics defined by (\ref{reduced-32}). 
It is easy to observe that $8\, f_{(32)}(x,y)=f_{(47)}(x,y)$ by a direct computation.

\vspace{.2cm}\noindent 
%\noindent 
{\bf nt67.} The moduli space $\cM (\Si_{67})$ with $\Si_{67}=[[2 A_2, A_5, E_6], [2 A_2]]$ 
is  not irreducible.
%  This is not obvious like the previous two cases. 
First we observe that $\emdim(\Si_{67})=0$ as before. 
We consider the minimal moduli slice with the slice condition: 
 
($\star$) two outer $A_2$'s are  at  $(0,\pm 1)$, an inner $E_6$ 
is at $(1,0)$  and an inner $A_5$ at $(-1,0)$. 
 
The corresponding slice reduces to two points defined by  
$\bff_a=(f_{2a},f_{3a})$ and $\bff_b=(f_{2b},f_{3b})$ where
\begin{eqnarray*}
\bff_a:\quad
\begin{cases}f_{2a}&=y^{2} + {  \frac {1}{2}}  - {  \frac {1}{
2}} \,x^{2} + {  \frac {1}{2}} \,I\,x^{2}\,\sqrt{3}
 - {  \frac {1}{2}} \,I\,\sqrt{3}\\
f_{3a}&={  \frac {1}{4}} \,\sqrt{18 - 6\,I\,\sqrt{3}}\,(1 - x
 + I\,\sqrt{3}\,y^{2} - x^{2} + x^{3} + x\,y^{2})
\end{cases}
\end{eqnarray*}
\begin{eqnarray*}
\bff_b:\quad \begin{cases}
f_{2b}&=y^{2} + {  \frac {1}{2}}  - {  \frac {1}{
2}} \,x^{2} - {  \frac {1}{2}} \,I\,x^{2}\,\sqrt{3}
 + {  \frac {1}{2}} \,I\,\sqrt{3}
\\
f_{3b}&={  \frac {1}{4}} \,\sqrt{18 + 6\,I\,\sqrt{3}}\,(1 - x
 - I\,\sqrt{3}\,y^{2} - x^{2} + x^{3} + x\,y^{2})
\end{cases}
\end{eqnarray*}
%Put $f^{(1)}:=f_{2a}^3+f_{3a}^2$ and $f^{(2)}:=f_{2b}^3+f_{3b}^2$.
%First, we claim that

\begin{Observation}   They are not in the same orbit of 
$\PGL(3,\bfC)$ in $\cM(\Si_{67})$.
\end{Observation}

{\em Proof.}
For  a matrix $A\in \GL(3,\bfC)$, we define as usual
$\phi_A: \bfP^2\to \bfP^2$ by the multiplication from the left.
Assume that  there is a matrix $A\in \GL(3,\bfC)$ such that 
$\bff_a^A:=\phi_A^*(\bff_a)=(\bff_b)$, it must keep the singular
points
$(-1,0),(1,0)$. Moreover  we observe that $f_{2a},f_{3a}, f_{2b},f_{3b}$
are even in $y$ variable. Thus the involution
$(x,y)\to (x,-y)$ keep the above polynomials. As the image of outer
singularities must be outer singularities,
 we may assume that $(0,1),(0,-1)$ are also  invariant by $\phi_A$.
This implies that $A=\Id$ in $\PGL(3,\bfC)$. This is ridiculous.\qed

\begin{Observation}\label{not-injective}  Each of $\psi_{red}(\bff_a),
\psi_{red}(\bff_b)$ has two different torus 
decompositions
in $\cM(\Si_{67})$.
\end{Observation}

{\em Proof.}  We will show the assertion for $\bff_a$.
First , two inner $A_2$ are located at
\[
P_1:= ({  \frac {-1}{3}} \,I\,\sqrt{3}, \,{  
\frac {1}{3}} \,\sqrt{3} + I),\qquad
P_2:=({  \frac {-1}{3}} \,I\,\sqrt{3}, \, - 
{  \frac {1}{3}} \,\sqrt{3} - I)
\]
We choose a conic
$h_2(x,y)=0$ which passes through four $A_2$ singularities
$(-1,0), (1,0), P_1,P_2$, and cut $x$-axis vertically at $(1,0)$.
Then another decomposition is given by
$\psi_{red}(\bff_a)=h_2^3+h_3^2$
where
$h_2(x,y):=
y^{2} - 1 + x^{2}$ and
\[h_3(x,y):=
{  \frac {1}{4}} \,(x^{3} - x^{2} - I\,y^{2}\,x\,
\sqrt{3} - x + 1 - y^{2})\,\sqrt{18 - 6\,I\,\sqrt{3}}
\]

\begin{Observation} \label{Observation}$\bfh=(h_2,h_3)$ and $\bff_b=(f_{2b},f_{3b})$ are in the same 
$\GL(3,\bfC)$-orbit in
$\tilde{\cM}(\Si_{67})$. In particular, $\psi_{red}(\bff_a)$ and
 $\psi_{red}(\bff_b)$ are $\PGL(3,\bfC)$-equivalent.
\end{Observation}
In fact,  a direct computation  shows that 
$\phi_B^*(h_2,h_3)=(f_{2b},f_{3b})$ where
\[B=
 \left[ 
{\begin{array}{ccc}
{  \frac {3}{4}}  & 0 & {  \frac {-1}{4}
} \,I\,\sqrt{3} \\ [2ex]
0 & {  \frac {1}{4}} \,\sqrt{3} + {  
\frac {3}{4}} \,I & 0 \\ [2ex]
{  \frac {-1}{4}} \,I\,\sqrt{3} & 0 & {  
\frac {3}{4}} 
\end{array}}
 \right] 
\]

\begin{Proposition}\label{Triple} The images  of the moduli spaces
$\cM([[4A_2,A_5],[E_6]])$, $\cM([[4A_2,E_6],[A_5]])$ and
 $\cM([[2A_2,A_5,E_6],[2A_2]])$  by the morphism
 $\psi_{red}$ into $\cM([4A_2,A_5,E_6])$
are the same. \end{Proposition}
%It is a direct computation to see that
The first equality $\psi_{red}(\cM([[4A_2,A_5],[E_6]]))=\psi_{red}(\cM([[4A_2,E_6],[A_5]]))$ 
is already observed  by the above normal forms.
Observation \ref{Observation} proves that $\cM([[2A_2,A_5,E_6],[2A_2]])$ is
irreducible.
Thus it is enough to show that 
$\psi_{red}(\bff_a)\in \psi_{red}(\cM([[4A_2,A_5],[E_6]])$.
In fact, we have 
$ \psi_{red}(\bff_a)=\psi_{red}(\bfg)$ where $\bfg=(g_2,g_3)\in
\psi_{red}(\cM([[4A_2,A_5],[E_6]])
$
 and $g_2,g_3$ are  given by

\begin{multline*}
\qquad{\mathrm{g_2}(x, \,y)=y^{2} - 1 + (1 + 2\,I\,\sqrt{3})\,x^{2} + 2
\,I\,\sqrt{3}\,x}
\\
{\mathrm{g_3}(x, \,y)={\displaystyle \frac {1}{28}} (7\,x\,y^{2}
 + 2\,y^{2} + I\,\sqrt{3}\,y^{2} + 4\,x^{3} + 9\,I\,x^{3}\,\sqrt{
3} - 2\,x^{2} + 13\,I\,x^{2}\,\sqrt{3} - 8\,x }\\
\mbox{} + 3\,I\,x\,\sqrt{3} - 2 - I\,\sqrt{3})\sqrt{ - 54 - 78\,I
\,\sqrt{3}} 
\end{multline*}

\vspace{.3cm} 
\noindent
{\bf {nt64.}}
We consider the distinguished configuration moduli
$\cM(\Si_{64})$ with $\Si_{64}=[[2A_2,A_5,E_6],[A_4]]$.
We have $\emdim(\Si_{64})=0$.
We consider the minimal slice with respect to:

($\star$) an inner  $A_5$ is at (0,1), an inner
 $E_6$ is at $(1,-1)$ with  tangent cone
$x=1$ and an outer $A_4$ is at $(0,0)$ with  tangent cone $y=0$.

The minimal slice consists of two points
$\bff_a=(f_{2a},f_{3a})$ and $\bff_b=(f_{2b},f_{3b})$:
\begin{eqnarray*}
\begin{cases}
&\mathrm{f_{2a}}(x, \,y)={\displaystyle \frac {1}{5}} \,
 (5\,y^{2}\,\sqrt{5} - 10 + 5\,y\,x + 4\,x\,
\sqrt{5} + 16\,x - y\,\sqrt{5} + 5\,y + 5\,x^{2}\,\sqrt{5} - x^{2
}\\
&\mbox{} - 4\,\sqrt{5} + 11\,y\,x\,\sqrt{5} + 5\,y^{2})/(1 + \sqrt{5}) 
\\
&\mathrm{f_{3a}}={\displaystyle \frac {1}{125}} \sqrt{50 + 30\,
\sqrt{5}}(250 + 110\,x^{3}\,\sqrt{5} + 88\,x^{3} - 420\,y\,x + 
110\,\sqrt{5} - 192\,y^{2}\,\sqrt{5} \\
&\mbox{} - 210\,x\,\sqrt{5} - 15\,y\,\sqrt{5} - 48\,x^{2}\,\sqrt{5
} + 155\,y^{3} + 366\,y\,x^{2} + 348\,y^{2}\,x\,\sqrt{5} + 336\,y
\,x^{2}\,\sqrt{5} \\
&\mbox{} + 97\,y^{3}\,\sqrt{5} - 510\,x - 75\,y + 498\,y^{2}\,x - 
300\,y\,x\,\sqrt{5} + 30\,x^{2} - 330\,y^{2}) \left/ {\vrule 
height0.45em width0em depth0.45em} \right. \!  \! (1 + \sqrt{5})
^{3}
\end{cases}
\end{eqnarray*}
\begin{eqnarray*}\begin{cases}
&\mathrm{f_{2b}}(x, \,y)={\displaystyle \frac {1}{5}} \,
(11\,y\,x\,\sqrt{5} - 5\,y\,x - y\,\sqrt{5}
 + 10 + 5\,y^{2}\,\sqrt{5} + x^{2} + 5\,x^{2}\,\sqrt{5} - 5\,y - 
16\,x \\&+ 4\,x\,\sqrt{5} - 4\,\sqrt{5} - 5\,y^{2})/( - 1 + \sqrt{5})\\
&\mathrm{f_{3b}}(x, \,y)={\displaystyle \frac {1}{125}} \,I\,
\sqrt{ - 50 + 30\,\sqrt{5}}( - 250 + 110\,x^{3}\,\sqrt{5} - 88\,x
^{3} + 420\,y\,x + 110\,\sqrt{5} \\
&\mbox{} - 192\,y^{2}\,\sqrt{5} - 210\,x\,\sqrt{5} - 15\,y\,\sqrt{
5} - 48\,x^{2}\,\sqrt{5} - 155\,y^{3} - 366\,y\,x^{2} + 348\,y^{2
}\,x\,\sqrt{5} \\
&\mbox{} + 336\,y\,x^{2}\,\sqrt{5} + 97\,y^{3}\,\sqrt{5} + 510\,x
 + 75\,y - 498\,y^{2}\,x - 300\,y\,x\,\sqrt{5} - 30\,x^{2} + 330
\,y^{2}) \\
 &\left/ {\vrule height0.48em width0em depth0.48em} \right. \! 
 \! ( - 1 + \sqrt{5})^{3} 
\end{cases}
\end{eqnarray*}

Note that the stabilizer in  $\PGL(3,\bfC)$
of  three points $(0,0)$, $(1,-1)$, $(0,1)$ and 
two lines $x=1$ and $y=0$ is trivial. Thus $\bff_a$ and $\bff_b$ are not 
in the same orbit even in the reduced moduli space
$\cM_{red}([2A_2,A_5,E_6,A_4])$. Thus the reduced moduli has two irreducible
components.
\begin{Proposition}\label{Galois1}
Two sextics $f_a:=f_{2a}^3+f_{3a}^2$ and
 $f_b:=f_{2b}^3+f_{3b}^2$
are defined over $\bfQ(\sqrt{5})$.
Let $\iota:\bfQ(\sqrt{5})\to\bfQ(\sqrt{5}) $ be the involution induced
 by the Galois automorphism defined by 
$\iota(\sqrt{5})=-\sqrt{5}$. Then  $\iota(f_a)=f_b$.
%we can easily see that 
%$f_a$ is mapped to $f_b$ by this involution.
\end{Proposition}
We do not know if  there exists an explicit homeomorphism of
the complements of 
 the sextics $f_a=0$ and $f_b=0$  in $\bfP^2$.
% or not.

\vspace{.3cm}
\noindent 
{\bf nt70.}  The moduli space $\cM(\Si_{70})$ with
$\Si_{70}=[[2 A_2, 2 E_6], [A_3]] $. The distinguished configuration moduli
is  irreducible and transversal
and $\emdim(\Si_{70})=\msdim(\Si_{70})=0$.
  For the computation of a minimal slice, we use the slice condition:

$(\star)$
an outer $A_3$ is  at the origin with tangent cone   $x=0$ and
 two inner $E_6$'s are at $(1,\pm 1)$. The tangent cone at $(1,1)$
is  given by   $y=1$.

The normal form is  given by 
\[
\begin{cases}
f_2(x,y)&=\frac{1}{3}(3y^2+(6x-6)y-2x^2-2x+1)\\ 
f_3(x,y)&=\frac{I\,\sqrt{3}}{9}(x-1)(18y^2+(9x-9)y-17x^2-2x+1)
\end{cases}\]

\vspace{.3cm}
\noindent
{\bf nt78.} We consider the moduli slice of $\cM(\Si_{78})$
where $\Si_{78}=[[3A_2,A_8],[D_5]]$.
We have $\emdim(\Si_{78})=0$. However the computation of minimal slice turns
out to be messy. So we consider the slice $\cA$ under the condition:

($\star$) an outer $D_5$ is at $O=(0,0)$ with  $y=0$  as the tangent cone of
multiplicity 1, and an inner
$A_8$ is  at $(1,1)$ with $y=1$ as the tangent cone.

The normal form $\bff=(f_2,f_3)$ is given by 
\begin{multline*}
\mathrm{f_2}(x, \,y)={\displaystyle \frac {1}{8\mathit{
t_1}^{2} }} (8\,\mathit{t_1}
^{4}\,y - 8\,\mathit{t_1}^{4} + 8\,y^{2}\,\mathit{t_1}^{2} + 8\,
\mathit{a_{10}}\,x\,\mathit{t_1}^{2} - 8\,y\,\mathit{a_{10}}\,x\,
\mathit{t_1}^{2} - 8\,y\,\mathit{t_1}^{2} + 2\,y\,\mathit{a_{10}}^{2}
\,x \\
\mbox{} - y\,\mathit{a_{10}}^{2} - \mathit{a_{10}}^{2}\,x^{2})
\end{multline*}
\begin{multline*}
\mathrm{f_3}(x, \,y)={\displaystyle \frac {1}{512}} ( - 24\,y^{3}
\,\mathit{a_{10}}^{2}\,\mathit{t_1}^{4} + 48\,\mathit{a_{10}}^{2}\,
\mathit{t_1}^{2}\,y^{2} + 192\,y^{2}\,\mathit{t_1}^{4} + 288\,
\mathit{t_1}^{4}\,x^{2}\,\mathit{a_{10}}^{2} + 512\,\mathit{t_1}^{8}
 \\
\mbox{} - 48\,y^{2}\,\mathit{a_{10}}^{2}\,x\,\mathit{t_1}^{2} - 16\,
\mathit{a_{10}}^{3}\,x^{3}\,\mathit{t_1}^{2} + 24\,y\,\mathit{a_{10}}^{3
}\,x^{2}\,\mathit{t_1}^{2} + 1152\,y\,\mathit{a_{10}}\,x\,\mathit{t_1}
^{6} \\
\mbox{} - 192\,y\,\mathit{a_{10}}^{2}\,x\,\mathit{t_1}^{4} - 48\,y\,
\mathit{a_{10}}^{3}\,x\,\mathit{t_1}^{2} + 3\,y\,\mathit{a_{10}}^{4}\,x
^{2} - 768\,y\,\mathit{t_1}^{8} + 24\,y\,\mathit{a_{10}}^{2}\,x^{2}\,
\mathit{t_1}^{2} \\
\mbox{} - 264\,y\,\mathit{t_1}^{4}\,x^{2}\,\mathit{a_{10}}^{2} + 384
\,y^{2}\,\mathit{a_{10}}\,x\,\mathit{t_1}^{4} + 48\,y^{2}\,\mathit{
a_{10}}^{3}\,x\,\mathit{t_1}^{2} - 384\,y^{2}\,\mathit{a_{10}}\,x\,
\mathit{t_1}^{6} \\
\mbox{} + 144\,y^{2}\,\mathit{a_{10}}^{2}\,x\,\mathit{t_1}^{4} - 48\,
\mathit{a_{10}}^{2}\,y^{2}\,\mathit{t_1}^{4} - 192\,y^{3}\,\mathit{t_1
}^{4} + 384\,\mathit{t_1}^{6}\,y^{3} - 768\,\mathit{a_{10}}\,x\,
\mathit{t_1}^{6} \\
\mbox{} + 64\,y^{3}\,\mathit{t_1}^{8} - 1152\,\mathit{t_1}^{6}\,y^{
2} + 768\,y\,\mathit{t_1}^{6} + 3\,\mathit{a_{10}}^{4}\,y^{2} - 8\,y
^{3}\,\mathit{t_1}^{2}\,\mathit{a_{10}}^{3} - 24\,y^{3}\,\mathit{t_1}
^{2}\,\mathit{a_{10}}^{2} \\
\mbox{} - 384\,y\,\mathit{a_{10}}\,x\,\mathit{t_1}^{4} + 192\,
\mathit{t_1}^{8}\,y^{2} - 6\,y^{2}\,\mathit{a_{10}}^{4}\,x + 96\,y\,
\mathit{t_1}^{4}\,\mathit{a_{10} }^{2})/\mathit{t_1}^{5} 
\end{multline*}
From this normal form, we
 see that $\cA$ is irreducible and we can fix one special 
point $\bff_a=(f_{2a},f_{3a})$, substituting 
$t_1=1,a_{10}=- 1$, where
\[
\mathrm{f_{2a}}(x, \,y)= - {\displaystyle \frac {1}{8}} \,y - 1 + y
^{2} - x + {\displaystyle \frac {5}{4}} \,x\,y - {\displaystyle 
\frac {1}{8}} \,x^{2}
\]
\[
\mathit{f_{3a}} := 1 - {\displaystyle \frac {57}{32}} \,x\,y - 
{\displaystyle \frac {261}{512}} \,x^{2}\,y + {\displaystyle 
\frac {21}{256}} \,x\,y^{2} + {\displaystyle \frac {1}{32}} \,x^{
3} - {\displaystyle \frac {765}{512}} \,y^{2} + {\displaystyle 
\frac {27}{64}} \,y^{3} + {\displaystyle \frac {3}{16}} \,y + 
{\displaystyle \frac {3}{2}} \,x + {\displaystyle \frac {9}{16}} 
\,x^{2}
\]

The isotropy subgroup fixing $(0,0),(1,1)$ and two lines
$y=0$ and $y=1$ is generated by
\[
A= \left[ 
{\begin{array}{ccc}
\mathit{a_1} & \mathit{a_2} & 0 \\
0 & \mathit{a_1} + \mathit{a_2} & 0 \\
0 & \mathit{a_1} + \mathit{a_2} - \mathit{a_9} & \mathit{a_9}
\end{array}}
 \right] 
\]
We can easily see that the orbit of $\bff_a$ by this isotropy group
is the whole slice $\cA$.
Thus $\cM(\Si_{78})$ has also a transversal minimal moduli slice
which is given by one point $\bff_a$.
In fact, we can see that ${\bf f}_a^A=\bf f$
where $A$ is defined by 
\[
 a_1= - {\displaystyle \frac {\mathit{a_{10}}}{\mathit{t_1}}},
\quad
a_2={\displaystyle \frac {1}{25}} \,{\displaystyle \frac {8\,\mathit{
t_1}^{4} + \mathit{a_{10}}^{2} + 17\,\mathit{a_{10}}\,\mathit{t_1}^{2} + 
8\,\mathit{t_1}^{2}}{\mathit{t_1}^{3}}},
\quad
a_9=\mathit{t_1}
\]
\vspace{.3cm}
\noindent 
{\bf nt83.} The moduli space $\cM(\Si_{83})$ with
$\Si_{83}=[[3 A_2, A_8], [A_1, A_4]] $ 
is  irreducible. Here we  compute the minimal slice $\cS$ with 
the following slice condition:

($\star$)
an outer $A_4$ is at the origin, an outer $A_1$ is at $(1,-1)$ 
and an inner $A_8$ is  at $(1,1)$. The  tangent cones at the origin
 and at $(1,1)$ are given by $x=0$ and  $y=1$ respectively.

Then $\emdim(\Si_{83})=0$ and 
 it has a transverse minimal slice which consists of  a single point
$\{(f_2,f_3)\}$
where 
\[
\begin{cases}
f_2(x,y)&={ \frac 1{565}}\,(565\,y^{2} + 126\,y\,x - 176\,y + 405\,x^{2}
 - 936\,x + 16)\\
f_3(x,y)&={  \frac {1}{319225}} \,I\,\sqrt{
565}(13321\,y^{3} + 28215\,y^{2}\,x - 6294\,y^{2} + 16767\,y\,x^{
2} - 31644\,y\,x \\
&\mbox{} + 1056\,y + 18225\,x^{3} - 45198\,x^{2} + 5616\,x - 64)
 \end{cases}
\]

\vspace{.3cm}
\noindent 
{\bf nt91.} We consider the moduli space $\cM(\Si_{91})$
with $\Si_{91}=[[3 A_2, B_{3,6}], [A_2]]$. 
The distinguished configuration moduli is  irreducible and transversal
and $\emdim(\Si_{83})=\msdim(\Si_{83})=2$.
For the computation of a minimal slice, we use the slice condition:

($\star$) an outer $A_2$ is  at $O=(0,0)$ with the tangent cone $x=0$, 
an inner $B_{3,6}$ is  at $(1,1)$ with the tangent cone $y=1$ and an
inner $A_2$ 
is at $(1,-1)$ 

The normal form are given by the following polynomials with 
 two-parameters $t_1,t_2$  
($t_1\ne 0$, $t_2\ne 0$):
\[
\begin{cases} 
f_2(x,y)&=y^2-(t_2x-t_2)y+(1+t_2-t_1^2)x^2+(2t_1^2-t_2-2)x-t_1^2\\ 
f_3(x,y)&=\frac{1}{8t_1}(6 t_2 y^3+((-6 t_2+12 t_1^2-3 t_2^2) x-12 t_1^2+3 t_2^2) y^2+((6 t_2+6 t_2^2-12 t_2 t_1^2) x^2\\
&\mbox{}+(-6 t_2^2+24 t_2 t_1^2-12 t_2) x-12 t_2 t_1^2) y+(-8 t_1^4+12 t_1^2-6 t_2-3 t_2^2+12 t_2 t_1^2) x^3\\
&\mbox{}+(3 t_2^2+24 t_1^4-24 t_2 t_1^2+12 t_2-36 t_1^2) x^2+(24 t_1^2-24 t_1^4+12 t_2 t_1^2) x+8 t_1^4)
\end{cases}
\]

\vspace{.3cm}
\noindent 
{\bf nt99.} The moduli space  $\cM (\Si_{99})$ with 
$\Si_{99}=[[A_5, 2 E_6], [A_2]]$ is not irreducible.
First we observe that $\emdim(\Si_{99})=0$ as before. 
We consider the minimal moduli slice with the slice condition:  

($\star$) an outer $A_2$ is  at $(-1,0)$, two inner  $E_6$'s are at $(0,\pm 1)$
 and
an inner $A_5$  is at
$(1,0)$.
 
The corresponding slice reduces to two points
% defined by  
$\bff_a=(f_{2a},f_{3a})$ and $\bff_b=(f_{2b},f_{3b})$ where
\[
\begin{cases}
f_{2a}(x,y)&= - {  \frac {1}{23}} \, (5 + 4\,\sqrt{3}) \,(5\,y^{2}
 - 5 + 23\,x^{2} - 18\,x - 4\,x\,\sqrt{3} - 4\,y^{2}\,\sqrt{3} + 
4\,\sqrt{3})\\
f_{3a}(x,y)&=2\,\sqrt{3 + 2\,\sqrt{3}}\,(1 + \sqrt{3})\,(\sqrt{3} + 3\,x^{2}
 - x\,\sqrt{3} - 3\,x - y^{2}\,\sqrt{3})\,x
\end{cases} 
\]
\[
\begin{cases}
f_{2b}(x,y)&={  \frac {1}{23}} \,( - 5 + 4\,\sqrt{3})\,(5\,y^{2}
 - 5 + 23\,x^{2} - 18\,x + 4\,x\,\sqrt{3} + 4\,y^{2}\,\sqrt{3} - 
4\,\sqrt{3})\\
f_{3b}(x,y)&= - 2\,\sqrt{3 - 2\,\sqrt{3}}\,( - 1 + \sqrt{3})\,( - \sqrt{3} + 3
\,x^{2} + x\,\sqrt{3} - 3\,x + y^{2}\,\sqrt{3})\,x
\end{cases} 
\] 
The isotropy subgroup fixing the configuration of singularity,
except possibly exchanging two $E_6$ is generated by the involution
$\iota\maketitle(x,y)\to (x,-y)$. However
 the defining conics and cubics are even in $y$. Thus $\bff_a,\bff_b$
are invariant under this involution.
Thus the moduli spaces $\cM(\Si_{99})$ and 
$\cM_{red}((\Si_{99})_{red})$ has  two irreducible components,
like the case nt64. Also we have a similar assertion:
\begin{Proposition}
Both sextics $\psi_{red}(\bff_a)=f_{2a}^3+f_{3a}^2=0$ and
 $\psi_{red}(\bff_b)=f_{2b}^3+f_{3b}^2=0$ are defined over $\bfQ(\sqrt{3})$.
Let $\iota:\bfQ(\sqrt{3})\to\bfQ(\sqrt{3}) $ be the involution induced
 by the Galois automorphism defined by 
$\iota(\sqrt{3})=-\sqrt{3}$. Then 
%we can easily see that 
$\iota(\psi_{red}(\bff_a))=\psi_{red}(\bff_b)$.
% by this involution.
\end{Proposition}

\vspace{.3cm}
\noindent 
{\bf nt100.} Let us consider the moduli space $\cM(\Si_{100})$ 
with  $\Si_{100}=[[3 E_6], [A_1]]$.
The distinguished configuration moduli is  irreducible and transversal
and $\emdim(\Si_{100})=\msdim(\Si_{100})=0$.
For the computation of a minimal slice, we use the slice condition:
  
($\star$) an outer $A_1$ is  at $(-1,0)$ and  three inner
$E_6$'s are at $(0,\pm 1)$, and $(1,0)$. 

The normal forms are given by
\[
%\begin{cases} 
f_2(x,y)=y^2-5x^2+6x-1,\quad
f_3(x,y)=6\sqrt{3}x(x+y-1)(x-y-1) 
\]

This curve has been studied in our previous  paper \cite{Oka-Pho1}.

\vspace{.3cm}
\noindent 
{\bf nt104.} The moduli space $\cM(\Si_{104})$ 
with  $\Si_{104}=[[A_2, A_5, A_8], [A_4]]$ is not irreducible.
First we observe that $\emdim(\Si_{104})=0$ as before. 
We consider the minimal moduli slice with the slice condition:  
 
($\star$) an outer $A_4$ is  at $(0,0)$ with the tangent cone $x=0$, 
an inner $A_8$ is at $(1,1)$
with the tangent cone $y=1$ and an inner $A_5$ is  at $(1,-1)$. 

The corresponding slice reduces to two points  %defined by  
$\bff_a=(f_{2a},f_{3a})$ and $\bff_b=(f_{2b},f_{3b})$ where
\[
\begin{cases}
{f_{2a}}(x,y) &= y^{2} + {  \frac {11}{5}} \,y
\,x - {  \frac {11}{5}} \,y - {  \frac {1
}{6}} \,x^{2} - {  \frac {28}{15}} \,x + 
{  \frac {31}{30}}  + I\,( - {  \frac {2
}{5}} \,y\,x + {  \frac {2}{5}} \,y - {  
\frac {1}{6}} \,x^{2} + {  \frac {11}{15}} \,x - 
{  \frac {17}{30}} )
\\
{f_{3a}}(x,y) &= {  \frac {1}{443682000}} 
\sqrt{-537594690 - 620415330\,I}( - 14148\,y^{3} + 7532 - 25008\,
x \\
&\mbox{} + 3925\,x^{3} - 41895\,y^{2}\,x - 21546\,y\,x^{2} + 72522
\,y\,x - 36828\,y + 12849\,x^{2} + 42597\,y^{2} \\
&\mbox{} + I( - 24093\,x^{2} - 18324 + 25497\,y\,x^{2} - 29529\,y
^{2} - 74754\,y\,x + 42696\,y + 43956\,x \\
&\mbox{} + 6561\,y^{3} + 27990\,y^{2}\,x)) 
\end{cases}\]
\[
 \begin{cases}
{f_{2b}}(x,y)&= y^{2} + {  \frac {11}{5}} \,y
\,x - {  \frac {11}{5}} \,y - {  \frac {1
}{6}} \,x^{2} - {  \frac {28}{15}} \,x + 
{  \frac {31}{30}}  - I\,( - {  \frac {2
}{5}} \,y\,x + {  \frac {2}{5}} \,y - {  
\frac {1}{6}} \,x^{2} + {  \frac {11}{15}} \,x - 
{  \frac {17}{30}} )\\
{f_{3b}}(x,y) &= {  \frac {1}{443682000}} 
\sqrt{-537594690 + 620415330\,I}( - 14148\,y^{3} + 7532 - 25008\,
x \\
&\mbox{} + 3925\,x^{3} - 41895\,y^{2}\,x - 21546\,y\,x^{2} + 72522
\,y\,x - 36828\,y + 12849\,x^{2} + 42597\,y^{2} \\
&\mbox{} - I( - 24093\,x^{2} - 18324 + 25497\,y\,x^{2} - 29529\,y
^{2} - 74754\,y\,x + 42696\,y + 43956\,x \\
&\mbox{} + 6561\,y^{3} + 27990\,y^{2}\,x)) 
\end{cases}
\]
\begin{Proposition}\label{Galois2}
Let $f_a:=f_{2a}^3+f_{3a}^2$ and $f_b:=f_{2b}^3+f_{3b}^2$ and
we consider the sextics $C_a:=\{f_a=0\}$ and  $C_b:=\{f_b=0\}$.
Let $\vphi : \bfC[x,y] \to \bfC[x,y]$
be the Galois involution defined by the complex
conjugation on the coefficients.
We first observe that $f_{b}=\vphi(f_a)$.
Let $\xi:\bfP^2\to \bfP^2$ be the  homeomorphism defined
 by  the complex conjugation $\xi((X,Y,Z)=(\bar X,\bar Y,\bar Z)$, or,
$\xi(x,y)=(\bar x,\bar y)$ in the affine coordinate.
The above observation gives the homeomorphism of the pairs of spaces
$\xi:(\bfP^2, C_a)\to (\bfP^2,C_b)$. In particular, their complements
$\bfP^2-C_a$ and $\bfP^2-C_b$ are homeomorphic.
\end{Proposition}

\vspace{.3cm}
\noindent 
{\bf nt110.} We consider the moduli space $\cM(\Si_{110})$ 
with  $\Si_{110}=[[A_2,E_6,A_8], [A_3]]$.
The distinguished configuration moduli is  irreducible and transversal
and $\emdim(\Si_{110})=\msdim(\Si_{110})=0$.
For the computation of a minimal slice, we use the slice condition:
  
($\star$) an outer $A_3$ is  at $(0,0)$ with the tangent cone $x=0$,  
an inner $A_8$ is at $(1,1)$ with the tangent cone $y=1$ and
an inner  $E_6$ is at $(1,-1)$. 

The defining polynomials are given by
%\[
\begin{eqnarray*} 
&\qquad f_2(x,y)=\frac{1}{15}(15\,y^{2} + 12\,y\,x - 12\,y + 5\,x^{2} - 22\,x + 2)\\   
\lefteqn{f_3(x,y)={  \frac {I\sqrt{30}}{450}} \, 
(81\,y^{3} + 180\,y^{2}\,x - 99\,y^{2} + 117\,y\,x^{2} - 234\,y\,
x + 36\,y + 40\,x^{3}} \\
&\mbox{}- 183\,x^{2} + 66\,x - 4)
\end{eqnarray*}

%\vspace{.3cm}
%\input nt113-normal.tex
\noindent 
{\bf nt113.} We consider the moduli space $\cM(\Si_{113})$
with $\Si_{113}=[[A_2, E_6, A_8], [A_1, A_2]]$.
First we observe that $\emdim(\Si_{113})=0$.
We first compute the minimal slice of $\cM([[2A_2,E_6,A_5],[A_1,A_1]])$ with respect to:

$(\star)$  an outer $A_1$ is at $P:=(0,-1)$, an outer $A_1$ is at $O:=(0,0)$,
an inner $E_6$ is at $Q:=(-1,-1)$ and an inner $A_5$ is at $R:=(-1,1)$.

The corresponding normal form is given by

\begin{eqnarray*}
%\begin{cases}
f_2(x,y)& =  - ( - y^{2}\,\mathit{t}\,\mathit{a_{01}} + 
\mathit{t}^{3}\,y^{2} + y^{2} + y^{2}\,\mathit{a_{01}} - y^{2}\,
\mathit{t} - \mathit{t}^{2}\,y^{2} - x\,\mathit{t}\,\mathit{
a_{01}}\,y + y\,x\,\mathit{a_{01}} \\
&\mbox{} - \mathit{a_{01}}\,y\,\mathit{t} + \mathit{a_{01}}\,y + x^{2}
\,\mathit{t}^{3} + 5\,x^{2}\,\mathit{t}^{2} - 2\,x^{2}\,
\mathit{t}\,\mathit{a_{01}} - 4\,x^{2}\,\mathit{a_{01}} - 6\,x^{2} - 3
\,x\,\mathit{t}\,\mathit{a_{01}} \\
&\mbox{} - 3\,\mathit{a_{01}}\,x - 6\,x + 4\,x\,\mathit{t}^{2} + 2\,
x\,\mathit{t}^{3} + \mathit{t} - 1) \left/ {\vrule 
height0.37em width0em depth0.37em} \right. \!  \! (\mathit{t} - 
1)\\
f_3(x,y)& =  - {  \frac {1}{2}} ( - 2 + 2
\,\mathit{t} + 12\,y\,x\,\mathit{t}^{2} + 6\,y^{2} + 4\,y^{3}
 - 20\,x^{2} - 3\,\mathit{a_{01}}\,y\,\mathit{t} - 9\,x\,\mathit{t
}\,\mathit{a_{01}} \\
&\mbox{} - 15\,x^{2}\,\mathit{t}\,\mathit{a_{01}} - 9\,y^{2}\,
\mathit{t}\,\mathit{a_{01}} + 9\,x^{2}\,\mathit{t}\,\mathit{a_{01}}\,
y - 8\,y\,x\,\mathit{t}^{3} + 16\,y\,x\,\mathit{t} - 9\,x^{2}\,
\mathit{t}^{2}\,\mathit{a_{01}} \\
&\mbox{} - 4\,y\,x\,\mathit{t}^{4} - 6\,y\,x\,\mathit{a_{01}} + 6\,
\mathit{t}^{2}\,y^{2} + 2\,\mathit{t}^{3}\,y^{3} - 16\,y\,x + 4
\,x^{2}\,\mathit{t}^{4} + 20\,x^{2}\,\mathit{t}^{3} + 6\,x^{2}
\,\mathit{t}^{2} \\
&\mbox{} - 20\,y\,x^{2} + 2\,\mathit{t}^{4}\,x^{3} + 16\,\mathit{
t}^{3}\,x^{3} - 2\,y^{2}\,x\,\mathit{t}^{3} - 2\,y^{2}\,x\,
\mathit{t} + 2\,y^{2}\,x\,\mathit{t}^{4} + 26\,y\,x^{2}\,
\mathit{t} \\
&\mbox{} - 2\,y\,x^{2}\,\mathit{t}^{4} - 10\,y\,x^{2}\,\mathit{t
}^{3} + 6\,y\,x^{2}\,\mathit{t}^{2} - 18\,\mathit{t}\,x^{3} + 2
\,y^{2}\,x - 10\,x^{2}\,\mathit{t} - 2\,\mathit{t}^{4}\,y^{3}
 \\
&\mbox{} + 6\,\mathit{t}^{2}\,y^{3} - 12\,\mathit{t}\,\mathit{
a_{01}}\,x^{3} - 12\,x^{2}\,\mathit{a_{01}}\,y + 3\,\mathit{a_{01}}\,y - 9
\,\mathit{a_{01}}\,x - 6\,y^{3}\,\mathit{t}\,\mathit{a_{01}} \\
&\mbox{} + 3\,y^{3}\,\mathit{t}^{2}\,\mathit{a_{01}} + 3\,y^{2}\,
\mathit{t}^{2}\,\mathit{a_{01}} + 3\,y^{2}\,\mathit{a_{01}}\,x - 12\,y
^{2}\,\mathit{t} + 6\,y^{2}\,\mathit{a_{01}} + 6\,x\,\mathit{t}^{3
} + 12\,x\,\mathit{t}^{2} \\
&\mbox{} - 12\,x^{2}\,\mathit{a_{01}} + 3\,y^{3}\,\mathit{a_{01}} - 10\,
y^{3}\,\mathit{t} + 3\,y\,x^{2}\,\mathit{t}^{2}\,\mathit{a_{01}}
 - 18\,x - 3\,y^{2}\,\mathit{t}\,\mathit{a_{01}}\,x \\
&\mbox{} - 6\,x^{3}\,\mathit{a_{01}}\,\mathit{t}^{2} + 6\,y\,x\,
\mathit{t}^{2}\,\mathit{a_{01}}) \left/ {\vrule 
height0.37em width0em depth0.37em} \right. \!  \! (\mathit{t} - 
1) 
%\end{cases}
\end{eqnarray*}
Now the conditions for  $R$ (respectively  $O$) to be $A_8$
(resp. $A_2$)
singularities are 
given by
\[
 A_8: \begin{cases}
&\mathit{g_1}=2\,\mathit{a_{01}} - 2\,\mathit{a_{01}}\,\mathit{s} - 
\mathit{a_{01}}\,\mathit{s}^{2} + 4 - 4\,\mathit{s} - 5\,\mathit{
s}^{2} + 3\,\mathit{s}^{3}=0\quad\text{or}
\\
&\mathit{g_2}=16 - 7\,\mathit{a_{01}}^{2} + 12\,\mathit{a_{01}} + 59\,
\mathit{s}^{5} - 49\,\mathit{s}^{4} - 23\,\mathit{s}^{6} + 3\,
\mathit{s}^{7} + 85\,\mathit{a_{01}}\,\mathit{s}^{2} + 2\,\mathit{
s}^{3} \\
&\mbox{} - 40\,\mathit{s}^{2} + 8\,\mathit{a_{01}}\,\mathit{s} + 56
\,\mathit{s} - 118\,\mathit{a_{01}}\,\mathit{s}^{3} + 44\,\mathit{
a_{01}}\,\mathit{s}^{4} - 2\,\mathit{a_{01}}\,\mathit{s}^{5} - 
\mathit{a_{01}}\,\mathit{s}^{6} \\
&\mbox{} + 23\,\mathit{a_{01}}^{2}\,\mathit{s} - 3\,\mathit{a_{01}}^{2}
\,\mathit{s}^{2} - 7\,\mathit{a_{01}}^{2}\,\mathit{s}^{3} + 2\,
\mathit{a_{01}}^{2}\,\mathit{s}^{4} =0
\end{cases}
\]
\begin{multline*}
A_2: \mathit{H_1}= - 20\,\mathit{s}^{6} + 120\,\mathit{s}^{5} + 12\,
\mathit{a_{01}}\,\mathit{s}^{5} - 12\,\mathit{a_{01}}\,\mathit{s}^{4}
 - 144\,\mathit{s}^{4} - 24\,\mathit{a_{01}}^{2}\,\mathit{s}^{3}
 - 448\,\mathit{s}^{3} \\
\mbox{} - 240\,\mathit{a_{01}}\,\mathit{s}^{3} + 768\,\mathit{a_{01}}
\,\mathit{s}^{2} + 1344\,\mathit{s}^{2} + 108\,\mathit{a_{01}}^{2}
\,\mathit{s}^{2} - 768\,\mathit{a_{01}}\,\mathit{s} - 1152\,
\mathit{s} \\
\mbox{} - 144\,\mathit{a_{01}}^{2}\,\mathit{s} + 192\,\mathit{a_{01}}
 + 48\,\mathit{a_{01}}^{2} + 256 =0
\end{multline*} where  we put $t=s-1$.
It turns out that $g_1=H_1=0$ gives three points, defined by
\begin{multline*}
\lefteqn{\mathrm{f_2}(x, \,y)=((x^{2} - 1 - 5\,y\,x + x - 5\,y + y^{2})\,
\mathit{s}^{2} + ( - 2 - 7\,x - 4\,y\,x - {\displaystyle \frac {
17}{2}} \,y^{2} - 4\,y + {\displaystyle \frac {7}{2}} \,x^{2})\,
\mathit{s}} \\
\mbox{} + 2 + {\displaystyle \frac {11}{2}} \,y^{2} - 2\,x^{2} + 
{\displaystyle \frac {11}{2}} \,y + {\displaystyle \frac {11}{2}
} \,y\,x + {\displaystyle \frac {11}{2}} \,x)/( - 2 + 2\,\mathit{
s} + \mathit{s}^{2})\end{multline*}
\begin{multline*}
\lefteqn{\mathrm{f_3}(x, \,y)=((6\,y\,x^{2} - 15\,y + 19\,x^{3} - 2 - 36
\,y^{2}\,x - 33\,y^{2} + 3\,y^{3} - 6\,y\,x + 3\,x + 21\,x^{2})\,
\mathit{s}^{2}} \\
\mbox{} + ( - 21\,x + 9\,y^{2}\,x - {\displaystyle \frac {33}{2}
} \,y^{2} - {\displaystyle \frac {51}{2}} \,y^{3} - 12\,y - 
{\displaystyle \frac {33}{2}} \,x^{2} - 4 - 25\,x^{3} - 21\,y\,x
 + {\displaystyle \frac {33}{2}} \,y\,x^{2})\,\mathit{s} \\
\mbox{} + 4 + {\displaystyle \frac {33}{2}} \,y^{3} + 30\,y^{2}
 + 3\,x^{2} - 12\,y\,x^{2} + {\displaystyle \frac {33}{2}} \,x + 
7\,x^{3} + {\displaystyle \frac {27}{2}} \,y^{2}\,x + 21\,y\,x + 
{\displaystyle \frac {33}{2}} \,y)/( \\
 - 4 + 4\,\mathit{s} + 2\,\mathit{s}^{2}) 
\end{multline*}
where $ - 1 + 2\,\mathit{s}^{3}=0$.
The other pair $g_2=H_1=0$  is equivalent to 
$g_2=0$ and 
$ 2\, s2^4-13\, s2^3+27\, s2^2-19\, s2+5=0$.
As $g_2$ has degree 2 in $a_{01}$, this gives 8 points.
Anyway we have that $\msdim(\Si_{113})=0$.
 
\vspace{.3cm}
\noindent
{\bf nt118.} We consider the moduli slice of $\cM(\Si_{118})$
where $\Si_{118}=[[2A_2,A_{11}],[A_4]]$.
We have $\emdim(\Si_{118})=0$. However the computation of minimal slice turns
out to be complicated. So we consider the slice $\cA$ under the condition:

($\star$) an inner  $A_{11}$ is at $O=(0,0)$ with the tangent cone $x=0$ and
an outer $A_4$ is at $(1,0)$ with $x=1$ as the tangent cone.

The normal form is given by 
\[
\begin{cases}
f_2(x,y)&=- {  \frac {864}{125}} \,
{  \frac {\mathit{a_{11}}^{2}\,y^{2}}{\mathit{a_{10}}^{4}}
}  + \mathit{a_{11}}\,x\,y + ( - \mathit{a_{10}} - {  
\frac {25}{576}} \,\mathit{a_{10}}^{4})\,x^{2} + \mathit{a_{10}}\,x\\
f_3(x,y)&={  \frac {1}{8640000}} ( - 
155271168\,y^{3}\,\mathit{a_{11}}^{3} + 59719680\,y^{2}\,\mathit{a_{10}
}\,\mathit{a_{11}}^{2}\,x \\
&\mbox{} + 34214400\,y^{2}\,\mathit{a_{10}}^{4}\,x\,\mathit{a_{11}}^{2}
 - 59719680\,y^{2}\,\mathit{a_{10}}\,\mathit{a_{11}}^{2} - 31104000\,
\mathit{a_{10}}^{5}\,x^{2}\,y\,\mathit{a_{11}} \\
&\mbox{} - 2700000\,y\,\mathit{a_{10}}^{8}\,x^{2}\,\mathit{a_{11}} + 
31104000\,\mathit{a_{10}}^{5}\,x\,y\,\mathit{a_{11}} + 2700000\,
\mathit{a_{10}}^{9}\,x^{3} \\
&\mbox{} + 8640000\,\mathit{a_{10}}^{6}\,x^{3} + 78125\,\mathit{a_{10}}
^{12}\,x^{3} - 2700000\,\mathit{a_{10}}^{9}\,x^{2} - 17280000\,x^{2}
\,\mathit{a_{10}}^{6} \\
&\mbox{} + 8640000\,\mathit{a_{10}}^{6}\,x)/\mathit{a_{10}}^{6}
\end{cases} 
\]
We can easily see that $\cA$ is irreducible and we can fix one special 
point $\bff_a=(f_{2a},f_{3a})$, substituting 
$a_{11}=a_{10}=1$, where
\[
\begin{cases}
f_{2a}(x,y)&= - {  \frac {864}{125}} \,y^{2}
 + y\,x - {  \frac {601}{576}} \,x^{2} + x\\
f_{3a}(x,y)&=- {  \frac {11232}{625}} \,y^{3
} + {  \frac {1359}{125}} \,y^{2}\,x - 
{  \frac {864}{125}} \,y^{2} - {  \frac {
313}{80}} \,y\,x^{2} + {  \frac {18}{5}} \,y\,x + 
{  \frac {18269}{13824}} \,x^{3} - {  
\frac {37}{16}} \,x^{2} + x
\end{cases}\]
The isotropy subgroup $G_0$ fixing $(0,0),(1,0)$ and two lines
$x=0$ and $x=1$ is  given by
\[
 G_0=
\left \{ \left[ 
{\begin{array}{crc}
u & 0 & 0 \\
0 & w & 0 \\
u - v & 0 & v
\end{array}}
 \right]\in \PGL(3,\bfC)
\,;\, \,u,v,w\in \bfC^*\right \}\]
We can also show that the orbit of $\bff_a$ by this isotropy group
is the whole slice $\cA$.
Thus $\cM(\Si_{118})$ has also a transversal minimal moduli slice
which is given by one point $\bff_a$.

\vspace{.3cm}
\noindent 
{\bf nt123.} For the  normal forms of $\cM(\Si_{123})$ with
$\Si_{123}=[[2A_2,C_{3,9}^\natural],[A_2]]$, see the next section.

\vspace{.3cm}
\noindent 
{\bf nt128.} Now we consider the moduli space $\cM(\Si_{128})$ 
with  $\Si_{128}=[[2 A_8], [A_3]]$.
The distinguished configuration moduli is  irreducible and transversal
and $\emdim(\Si_{128})=0$.
For the computation of a minimal slice, we use the slice condition:
  
($\star$) an outer $A_3$ is  at $(-1,0)$ with the tangent cone $x=-1$, 
an inner $A_8$ is at $(0,1)$ with the tangent cone $y=1$ and another
inner  $A_8$ is at $(0,-1)$. 

The defining polynomials are given by
\[
\begin{cases} 
f_2(x,y)&=-3y^2-6xy-x^2+6x+3\\   
f_3(x,y)&= \frac{1}{16}(81 y^3+252 y^2 x+207 x^2 y-162 x y-81 y+38 x^3-180 x^2-90 x)
\end{cases}
\]

\vspace{.3cm}
\noindent 
{\bf nt136.} We consider the moduli space $\cM(\Si_{136})$ 
with  $\Si_{136}=[[E_6, A_{11}], [A_2]]$.
The distinguished configuration moduli is  irreducible and transversal
and $\emdim(\Si_{136})=\msdim(\Si_{136})=0$.
For the computation of a minimal slice, we use the slice condition:

($\star$) an inner  $A_{11}$ at $(0,0)$ with the tangent cone $x=0$, 
an outer $A_2$ at $(1,1)$ with the tangent cone $y+x=2$ and an inner 
$E_6$ at $(1,-1)$. 

The normal form is given by 

\begin{eqnarray*}
%\begin{multline*}
& \mathrm{f_2}(x, \,y)=y^{2} + {\displaystyle \frac {4}{3}} \,x\,y
 - {\displaystyle \frac {11}{3}} \,x^{2} + 4\,x\qquad\qquad\\
&{ \mathrm{f_3}(x, \,y)={\displaystyle \frac {1}{36}} \,I\,(14\,y^{3
} + 18\,y^{2}\,x + 12\,y^{2} - 54\,x^{2}\,y + 72\,x\,y - 10\,x^{3
} - 36\,x^{2} + 48\,x)\,\sqrt{6}}
\end{eqnarray*}

%\vspace{.3cm} 

%\input nt139-normal.tex
\noindent
{\bf nt139.} We consider the moduli slice of $\cM(\Si_{139})$
where $\Si_{139}=[[A_2,A_{14}],[A_3]]$.
We have $\emdim(\Si_{139})=0$. 
We consider the slice $\cA$ under the condition:

($\star$) an inner $A_{14}$ is at $O=(0,0)$ with the tangent cone $x=0$ and
an outer $A_3$ at $(1,0)$ with $x=1$ as the tangent cone, and an inner $A_2$
at $(-1,-1)$.

The corresponding slice is reduced to a single point
and  we can take
the normal form as follows.
\begin{eqnarray*}
&\mathrm{f_2}(x, \,y)=y^{2} - {\displaystyle \frac {10}{3}} \,x\,y
 + {\displaystyle \frac {41}{18}} \,x^{2} - {\displaystyle 
\frac {1}{18}} \,x
\\
&\mathrm{f_3}(x, \,y)= - {\displaystyle \frac {7}{16}} \,I\,\sqrt{
5}\,y^{3} + {\displaystyle \frac {433}{192}} \,I\,y^{2}\,\sqrt{5}
\,x - {\displaystyle \frac {1}{192}} \,I\,y^{2}\,\sqrt{5} - 
{\displaystyle \frac {27}{8}} \,I\,y\,\sqrt{5}\,x^{2} + 
{\displaystyle \frac {1}{24}} \,I\,y\,\sqrt{5}\,x \\
&\mbox{} + {\displaystyle \frac {1771}{1152}} \,I\,\sqrt{5}\,x^{3}
 - {\displaystyle \frac {97}{1728}} \,I\,\sqrt{5}\,x^{2} + 
{\displaystyle \frac {1}{3456}} \,I\,\sqrt{5}\,x 
\end{eqnarray*}

\vspace{.3cm} 
\noindent
{\bf nt142.} We consider the moduli slice of $\cM(\Si_{142})$
where $\Si_{142}=[[A_2,A_{14}],[A_1,A_2]]$.
We have $\emdim(\Si_{142})=0$. The  minimal slice
under the condition:

($\star$) an inner $A_{14}$ is at $O=(0,0)$ with the tangent cone $x=0$, an
outer $A_2$ is at $(1,0)$ with the  tangent cone
$x=1$ and an outer $A_1$ at $(-1,1)$.

The normal form is given by one point described by
\begin{eqnarray*}
&\mathrm{f_2}(x, \,y)=y^{2} + {\displaystyle \frac {16}{3}} \,x\,
y + {\displaystyle \frac {106}{45}} \,x^{2} - 2\,x
\\
&\mathit{f_3} = {\displaystyle \frac {41}{27}} \,I\,y^{3}\,
\sqrt{5} + {\displaystyle \frac {403}{45}} \,I\,y^{2}\,x\,\sqrt{5
} - {\displaystyle \frac {5}{9}} \,I\,y^{2}\,\sqrt{5} + 
{\displaystyle \frac {122}{15}} \,I\,y\,x^{2}\,\sqrt{5} - 6\,I\,y
\,x\,\sqrt{5} + {\displaystyle \frac {1354}{675}} \,I\,x^{3}\,
\sqrt{5} \\
&\mbox{} - {\displaystyle \frac {136}{45}} \,I\,x^{2}\,\sqrt{5} + 
{\displaystyle \frac {10}{9}} \,I\,\sqrt{5}\,x 
\end{eqnarray*}

\vspace{.3cm} 
\noindent
{\bf nt145.} We consider the moduli slice of $\cM(\Si_{145})$
where $\Si_{145}=[[A_{17}],[A_2]]$.
We have $\emdim(\Si_{145})=0$. However the computation of minimal slice turns
out to be complicated. So we consider the slice $\cA$ under the condition:

($\star$) $A_{17}$ is at $O=(0,0)$ with the tangent cone $x=0$ and
$A_2$ is at $(1,0)$.

We note  that the tangent cone at $A_2$ can not be generic.
In fact,  we see,  by computation, that the tangent cone at $A_2$ must
pass through  $A_{17}$.
The normal form is given by three dimensional family:
%$a_{10},b_{02},b_{11}$ ($a_{10}\ne 0,b_{02}\ne 0$):

\[
\begin{cases}
f_2(x,y)&=\mathit{a_{10}}\mathit{b_{02}}\,y^{2} + 
{  \frac {1}{2}} \,\mathit{a_{10}}\,\mathit{b_{11}}\,x\,y
 + ( - \mathit{a_{10}} - {  \frac {9}{64}} \,\mathit{a_{10}
}^{4})\,x^{2} + \mathit{a_{10}}\,x\\
f_3(x,y)&={  \frac {1}{2}}  \,\mathit{b_{02}}
\mathit{b_{11}}\,y^{3} - {  \frac {27}{64}} \,y^{2}\,x\,
\mathit{a_{10}}^{3}\mathit{b_{02}}  -\mathit{b_{02}}\, y^{2}\,x  + 
{  \frac {1}{4}} \,y^{2}\,x\,\mathit{b_{11}}^{2} + 
\,\mathit{b_{02}}\,y^{2} \\
&\mbox{} - {  \frac {9}{32}} \,\mathit{b_{11}}\,x^{2}\,y
\,\mathit{a_{10}}^{3} - \mathit{b_{11}}\,x^{2}\,y + \mathit{b_{11}}\,x\,y
 + x^{3} + {  \frac {9}{16}} \,x^{3}\,\mathit{a_{10}}^{3
} + {  \frac {27}{512}} \,x^{3}\,\mathit{a_{10}}^{6} \\
&\mbox{} - {  \frac {9}{16}} \,x^{2}\,\mathit{a_{10}}^{3}
 - 2\,x^{2} + x 
\end{cases} 
\]
We can easily see that $\cA$ is irreducible and we can fix one special 
point $\bff_a=(f_{2a},f_{3a})$, substituting $b_{11}=0$ and
$a_{10}=b_{02}=1$, where
\[
%\begin{cases}
f_{2a}(x,y)= y^{2} - {  \frac {73}{64}} \,x^{2} + x,\quad
f_{3a}(x,y)= - {  \frac {91}{64}} \,x\,y^{2}
 + y^{2} + {  \frac {827}{512}} \,x^{3} - 
{  \frac {41}{16}} \,x^{2} + x
%\end{cases}
\]
The isotropy subgroup $J$ fixing $(0,0),(1,0)$ and one lines
$x=0$ is 3-dimensional and it is  given  by
\[J=\left\{
M=\left[ 
{\begin{array}{ccr}
v + s & 0 & 0 \\
0 & u & 0 \\
v & w & s
\end{array}}
\right] \in \PGL(3,\bfC);\, u,\,s\ne 0,v\ne -1\right\}
\]
We can also show that the orbit of $\bff_a$ by this isotropy group
is the whole slice $\cA$.
Thus $\cM(\Si_{145})$ has also a transversal minimal moduli slice
which is given by one point $\bff_a$.

\section{Coincidence of some moduli spaces}\label{SameModuli}
We have seen that there exist $121(=145-24)$ different distinguished
configurations.
On the other hand,  we assert 
\begin{Theorem}\label{same-moduli}
For the following six reduced configurations, the corresponding distinguished
 configurations are not unique:
$[6A_2,A_5],\,  [6A_2,E_6]$,
 $[6A_2,A_1,A_5]$, $[4A_2,2 A_5],\, [4A_2,A_5,E_6],\, [3A_2,C_{3,9}]$.
 More precisely, we have
\begin{enumerate}
\item
   $\psi_{red}([[6A_2],[A_5]])=
\psi_{red}([[4A_2,A_5],[2A_2]])$  ( nt5 and nt37).
%where nt5: $[[6A_2],A5]$ and nt37: $[[4A_2,A_5],[2A_2]]$.
\item  $\psi_{red}([[6A_2],[E_6]])=\psi_{red}([[4A_2,E_6],[2A_2]])$ (nt8, nt52).
%nt8: $[[6A_2],[E_6]]$ and nt52: $[[4A_2,E_6],[2A_2]]$.
\item   $\psi_{red}([[6A_2], A_1,A_5])=\psi_{red}([[4A_2,A_5],A_1,2A_2])$ (nt 13, nt42).
%nt13:$[[6A_2], A_1,A_5]$ and nt42: $[[4A_2,A_5],A_1,2A_2]$.
\item  $\psi_{red}([[4A_2,A_5],[A_5]])$
$=\psi_{red}([[2A_2,2A_5], 2A_2])$ (nt29,  nt60).
%nt29: $[[4A_2,A_5],[A_5]]$ and nt60: $[[2A_2,2A_5], 2A_2]$.
\item   $\psi_{red}([[4A_2,A_5],E_5])=\psi_{red}([[4A_2,E_6],[A_5]])
=\psi_{red}([[2A_2,A_5,E_6],[2A_2]])$ (nt32, nt47, nt67).
%nt32:$[[4A_2,A_5],E_5]$, nt47: $[[4A_2,E_6],[A_5]]$ and nt67:
%$[[2A_2,A_5,E_6],[2A_2]]$.
\item  $\psi_{red}([[2A_2,C_{3,9}^\natural],[A_2]])
=\psi_{red}([3A_2,C_{3,9}])$  ( nt123 and t11 ).
\end{enumerate}
\end{Theorem}

{\em Proof.} We prove the assertion by giving
 explicit torus decompositions for a given $f\in \cM_{red}(\Si)$ using 
minimal moduli slices.

\vspace{.3cm}\noindent
{\bf I}. We will show that the 
respective images of  $\cM(\Si_5)$ and $\cM(\Si_{37})$ into
  the  reduced moduli space $\cM_{red}([6A_2,A_5])$
 coincide, where
$\Si_5=[[6A_2],[A_5]])$ and $\Si_{37}=[[4A_2,A_5],[2A_2]])$.
As their minimal slice dimensions  are both equal to  two, 
this case requires a 
heavy  computation. 
So we need a special device for the computation.
We first compute the normal form of the minimal moduli slice of 
$\cM(\Si_5)$, with the slice conditions:

($\star_1$): an outer $A_5$ is at $O:=(0,0)$  with the tangent cone $x=0$.

($\star_2$) Two inner $A_2$ are at $P:=(1,1)$ and $ Q:=(1,-1)$.
The tangent cone at $P$ is given by $y=1$.

First, we can easily observe that $\cM(\Si_5)$ is irreducible, by
looking at the slice with respect to $(\star_1)$.
Then we compute the minimal slice with respect to $(\star_1+\star_2)$.
There are  several components but we can use the following component
$\cA$
by the irreducibility of $\cM(\Si_5)$.
\begin{eqnarray}\label{nt5-normal}
 \begin{cases}
&\mathrm{f_2}(x, \,y)=y^{2} + ( - 1 - \mathit{a_{10}} + \mathit{t_0}^{
2})\,x^{2} + \mathit{a_{10}}\,x - \mathit{t_0}^{2}
\\
&\mathrm{f_3}(x, \,y)= - {\displaystyle \frac {1}{2}} ( - 3\,y^{2}
\,x\,\mathit{t_0}^{2} + 3\,y^{2}\,\mathit{a_{10}}\,x + 6\,y^{2}\,x - 
3\,\mathit{t_0}^{2}\,y^{2} + 4\,x^{3}\,\mathit{t_0}^{4} - 9\,x^{3}
\,\mathit{t_0}^{2}\\
& - 3\,x^{3}\,\mathit{a_{10}}\,\mathit{t_0}^{2} 
 + 3\,x^{3}\,\mathit{a_{10}} + 6\,x^{3} - 6\,x^{2}\,\mathit{
t_0}^{4} + 15\,x^{2}\,\mathit{t_0}^{2} + 6\,x^{2}\,\mathit{a_{10}}\,
\mathit{t_0}^{2} - 6\,x^{2}\,\mathit{a_{10}}\\ 
&\qquad\qquad - 12\,x^{2} 
 - 3\,\mathit{a_{10}}\,\mathit{t_0}^{2}\,x + 2\,\mathit{t_0}^{4
})/\mathit{t_0}
\end{cases}
\end{eqnarray}
Note that $t_0=f_3(0,0)/f_2(0,0)$.
We observe that $f_2(x,y), f_3(x,y)$ are even in $y$-variable
and $t_0$ is even in $f_2(x,y)$ and  in $t_0\,f_3(x,y)$.
Thus 
 the sextics $f_2^3+f_3^2=0$ is  symmetric with respect
to
$x$-axis and the change $t_0\to -t_0$ does not change the class
of $(f_2,f_3)$ in
$\cM(\Si_5)$.
In fact, this is the reason we consider the above slice condition.
For the computation of the minimal slice $\cM(\Si_{37})$, we consider
the slice $\cB$ with the  condition:

$(\star_3)$  Two outer $A_2$ at $P,Q$ and an inner $A_5$ at $O$. The
tangent cone at $O$ and $P$ are given by $x=0$ and $y=1$.

The normal form is given by $g(x,y)=g_2(x,y)^3+g_3(x,y)^2$ where
\begin{eqnarray}\label{nt37-normal}
\begin{cases}
&\mathrm{g_2}(x, \,y)=y^{2} + \mathit{a_{20}}\,x^{2} + ( - 1 - 
\mathit{a_{20}} - \mathit{t_1}^{2})\,x
\\
&\mathrm{g_3}(x, \,y)= - {\displaystyle \frac {1}{8}} \,
%{\displaystyle \frac {
( - 6\,x\,\mathit{t_1}^{2} + 6\,\mathit{a_{20}}
\,x - 6\,\mathit{a_{20}} + 6 - 6\,x - 6\,\mathit{t_1}^{2})\,y^{2}/
\mathit{t_1}\\
&  - {\displaystyle \frac {1}{8}} (6\,x^{2}\,\mathit{
t_1}^{2} - 6\,\mathit{a_{20}}\,x^{2} 
 + 3\,x\,\mathit{t_1}^{4} + 6\,x\,\mathit{t_1}^{2} + 6\,
\mathit{a_{20}}\,x - 9\,x - 3\,x^{3} + 3\,x^{3}\,\mathit{a_{20}}^{2} \\
&- 6\,x^{3}\,\mathit{a_{20}}\,\mathit{t_1}^{2} - x^{3}\,\mathit{t_1}^{4}
+ 6\,x^{2}\,\mathit{t_1}^{4} + 12\,x^{2} - 6\,x^{2}\,
\mathit{a_{20}}^{2} + 3\,x\,\mathit{a_{20}}^{2} + 6\,x\,\mathit{a_{20}}\,
\mathit{t_1}^{2})/\mathit{t_1} 
\end{cases}
\end{eqnarray}
Here $t_1=f_3(P)/f_2(P)$.
We observe that $g_2,g_3$ are also even in $y$-variable, while $t_1$ is
even in $f_2$ and  in $t_1\,f_3$.
The assertion follows from
\begin{Proposition}\label{morphisms} There are canonical 
bijective morphisms
$\xi_1: \cA\to \cB$ and $\xi_2:\cB\to \cA$ so that $\xi_1\circ \xi_2$
and $\xi_2\circ \xi_1$ induce  the identity maps on the images $\pi(\cA)$ and $\pi(\cB)$.
\end{Proposition}
{\em Proof.} First we construct $\xi_1$.
Take a  $\bff_a=(f_2,f_3)$ in $\cA$ written as (\ref{nt5-normal}). First we
show the 
existence
of a conic $g_2(x,y)=0$ which contains 4 $A_2$ singularities
of $f_2^3+f_3^2=0$ other than $P,Q$
and $A_5$ with the tangent line  $y=0$ at $O$. 
Four $A_2$ are symmetric with respect to $x$-axis and their
x-coordinates are the solutions of 
\[
 R_1=3\,{x}^{2}{{\it t_0}}^{2}+6\,{\it b_{12}}\,{x}^{2}{\it t_0}+4\,{x}^{2}{{
\it b_{12}}}^{2}+3\,x{{\it t_0}}^{2}+6\,{\it b_{12}}\,x{\it t_0}+3\,{{\it t_0}}
^{2}=0
\]
We do not need  to solve these solutions explicitly.
We start from the form $h_2(x,y)=y^2+a\, x^2+b\, x+c$. First we put the
condition
$h_2(0,0)=0$. Then we compute the resultant $S(x)$ of $h_2$ and $f_3$ in $y$.
Then by the above symmetry condition, $S$ can be written as 
$S(x)=S_1(x)^2$ where $S_1$ is a polynomial of degree 3.
%$S(x)=x^2$ and the other solutions
%are  the roots of  a quadratic polynomial $S_1(x)$ of $x$.
Then  $S_1$ must be  divisible by $R_1$.
This  condition is
enough to solve the coefficient of $h_2$ up to a multiplication of a
constant,
 and we have
\begin{multline*}
\mathrm{h_2}(x, \,y)=(4\,\mathit{t_0}^{4}\,x + 8\,x^{2}\,\mathit{
t_0}^{4} - 19\,x^{2}\,\mathit{t_0}^{2} - 2\,\mathit{a_{10}}\,\mathit{
t_0}^{2}\,x + \mathit{t_0}^{2}\,y^{2} - 10\,x^{2}\,\mathit{a_{10}}\,
\mathit{t_0}^{2} - 6\,\mathit{t_0}^{2}\,x \\
\mbox{} + 12\,x^{2}\,\mathit{a_{10}} + 12\,x^{2} + 3\,\mathit{a_{10}}^{
2}\,x^{2})/\mathit{t_0}^{2} 
\end{multline*}
Now we have to find the partner cubic polynomial $g_3(x,y)$
such that $f(x,y)=h_2(x,y)^3\,t+h_3(x,y)^2$ for some polynomial $h_3(x,y)$.
The argument  by Tokunaga \cite{Tokunaga-torus} can not be used as we
have $A_5$. Instead of using that, we introduce a systematic
computational method.
For that purpose, we consider the family of polynomial
$f_t(x,y):=f(x,y)-h_2(x,y)^3\,t$.  Assuming the existence of such $h_3$,
this family of sextics $f_t=0$ has 4 $A_2$ at the same location as $f=0$
and $A_5$ at the origin. (Note that the tangent line of the conic
$h_2=0$
at $O$ is the same with that of $f_3=0$.)
If there is a $\tau_0$ such that
$f_{\tau_0}$ is a square of a cubic polynomial, $f_{\tau_0}(x,y)=0$ has an
non-isolated singularity at $O$. 
So we look for a special value for which the singularity
at $O$ is bigger than $A_5$.  In fact  such a $\tau_0$ is given by
$\tau_0=1$ and then  we see that 
$f_{\tau_0}$ is a square of 
a polynomial of degree 3.
This technique is quite useful to find the partner cubic for other cases
and hereafter we refer this technique as {\em degeneration method}.
The corresponding cubic form is given by
\begin{multline*}
\mathrm{h_3}(x, \,y)={\displaystyle \frac {-1}{2}} \,I(\mathit{t_0
}^{4}\,y^{2} - 4\,\mathit{t_0}^{4}\,x + 2\,\mathit{t_0}^{6}\,x - 6
\,y^{2}\,x\,\mathit{t_0}^{2} + 5\,y^{2}\,x\,\mathit{t_0}^{4} + 36\,
x^{2}\,\mathit{t_0}^{2} - 53\,x^{2}\,\mathit{t_0}^{4} \\
\mbox{} + 20\,x^{2}\,\mathit{t_0}^{6} - 48\,x^{3} + 114\,x^{3}\,
\mathit{t_0}^{2} - 93\,x^{3}\,\mathit{t_0}^{4} + 26\,x^{3}\,
\mathit{t_0}^{6} - \mathit{t_0}^{4}\,x\,\mathit{a_{10}} - 3\,y^{2}\,
\mathit{a_{10}}\,\mathit{t_0}^{2}\,x \\
\mbox{} + 30\,x^{2}\,\mathit{a_{10}}\,\mathit{t_0}^{2} - 22\,x^{2}\,
\mathit{t_0}^{4}\,\mathit{a_{10}} - 72\,x^{3}\,\mathit{a_{10}} + 117\,x
^{3}\,\mathit{a_{10}}\,\mathit{t_0}^{2} - 49\,x^{3}\,\mathit{t_0}^{4}
\,\mathit{a_{10}} \\
\mbox{} + 6\,x^{2}\,\mathit{t_0}^{2}\,\mathit{a_{10}}^{2} - 36\,x^{3}
\,\mathit{a_{10}}^{2} + 30\,x^{3}\,\mathit{a_{10}}^{2}\,\mathit{t_0}^{2}
 - 6\,\mathit{a_{10}}^{3}\,x^{3})\sqrt{3}/\mathit{t_0}^{3}
\end{multline*}
Thus we define $\xi_1(f_2,f_3)=(h_2,h_3)$.
In terms of the parameters, $\xi_1$ is defined by 
$ \xi_1 (a_{10},t_0)= (a_{20},t_1)$ where
\begin{multline*}
\mathit{a_{20}} = {\displaystyle \frac {3\,\mathit{a_{10}}^{2} - 10
\,\mathit{a_{10}}\,\mathit{t_0}^{2} + 12\,\mathit{a_{10}} - 19\,\mathit{
t_0}^{2} + 8\,\mathit{t_0}^{4} + 12}{\mathit{t_0}^{2}}},\quad
\mathit{t_1} = {\displaystyle \frac {I\,\sqrt{3}\,(\mathit{a_{10}
} - 2\,\mathit{t_0}^{2} + 2)}{\mathit{t_0}}}
\end{multline*}

Now the construction of the morphism $\xi_2:\cB\to \cA$ is done in exact
same way. Take
$\bfg=(g_2,g_3)\in \cB$ as in (\ref{nt37-normal}). First find a conic which pass through 
$6A_2$ of $g(x,y)=0$, and then find the partner cubic by degeneration
method.
In term of parameters, we define $\xi_2(a_{20},t_1)=(a_{10},t_0)$ where
\[
\mathit{a_{10}} =  - {\displaystyle \frac {1}{2}} \,
\displaystyle \frac {3\,\mathit{a_{20}}^{2} + 5\,\mathit{a_{20}}\,
\mathit{t_1}^{2} - 6\,\mathit{a_{20}} - \mathit{t_1}^{2} + 3 + 2\,
\mathit{t_1}^{4}}{\mathit{t_1}^{2}},\quad
\mathit{t_0} = {\displaystyle \frac {-1}{2}}\displaystyle \frac { \,I\,\sqrt{3}\,( - 1 + \mathit{t_1}^{2} + \mathit{a_{20}})}{
\mathit{t_1}}
\]

We can easily check that 
  $\xi_1\circ\xi_2(g_2,g_3)=(g_2,-g_3)$ and
$\xi_2\circ\xi_1(f_2,f_3)=(f_2,-f_3)$ which implies the assertion.
(Recall that $(f_2,f_3)\sim (f_2,-f_3)$.)
\qed
\begin{Remark}
We remark that the generic element of $\cA$  is contained in the moduli
 space
$\cM(\Si_{5})$. However for non-generic element
$(f_2,f_3)\in \cA$, the slice condition $(\star_1+\star_2)$ guarantee
only  that $f_2^3+f_3^2=0$  has an outer $A_5$ at $O$ and two inner
 $A_2$ at $P,Q$.
As $ \Si_{13}$ or $\Si_{29}$ has 
 an outer $A_5$ and 4 inner $A_2$, their slices
with respect to the slice condition $(\star_1+\star_2)$ are subvarieties
 of $\cA$. Here $\Si_j$ is the configuration corresponding to nt-j in
 the table at the end.
Similarly the slices of $\cM(\Si_{42}),~\cM(\Si_{60})$ with respect to
 the slice condition $(\star_3)$ are subvarieties of $\cB$.
This observation will be used in the next two pairs.
\end{Remark}

%\input SameModuliIII.tex 
 
%\vspace{.3cm}the
\noindent 
{\bf II}. The equalities  
$ \psi_{red}(\cM(\Si_{13}))=\psi_{red}(\cM(\Si_{42}))$ and
$ \psi_{red}(\cM(\Si_{29}))=\psi_{red}(\cM(\Si_{60}))
$
follow from the above argument ( Proposition \ref{morphisms}),
where $\Si_{13}=[[6A_2],[A_1,A_5]]$,
$\Si_{42}=$\newline $[[4A_2,A_5],[A_1,2A_2]]$, 
$\Si_{29}=[[4A_2,A_5],[A_5]]$ and $\Si_{60}=[[2A_2,2A_5],[2A_2]]$.
First we consider the equality
 $\psi_{red}(\cM(\Si_{13}))=\psi_{red}(\cM(\Si_{42}))$.
In fact, we may consider the slice  $\cA', \cB'$ of 
$\cM(\Si_{13})$ or $\cM(\Si_{42})$ 
subject to the slice condition $(\star_1+\star_2)$ or $(\star_3)$.
Then we have the canonical inclusions
$ \cA'\subset \cA,\quad \cB'\subset \cB$.
For example, $\cA'$ consist of $(f_2,f_3)\in \cA$ such that 
$f_2^3+f_3^2=0$ has also an outer $A_1$. As $f_2,f_3$ are symmetric with
respect
to $x$-axis, $A_1$ must be on $y=0$.  Thus the condition for
$(f_2,f_3)$ to be in $ \cA'$ is given by the vanishing of the discriminant
polynomial of $f(x,0)/x^2$ in $x$, which is 
\begin{multline*}
( - 8 + 12\,\mathit{t_0}^{2} + 3\,\mathit{a_{10}
}\,\mathit{t_0}^{2} - 6\,\mathit{t_0}^{4})(3\,\mathit{a_{10}}^{4} - 24
\,\mathit{a_{10}}^{3}\,\mathit{t_0}^{2} + 30\,\mathit{a_{10}}^{3} \\
\mbox{} + 72\,\mathit{t_0}^{4}\,\mathit{a_{10}}^{2} - 180\,\mathit{t_0
}^{2}\,\mathit{a_{10}}^{2} + 120\,\mathit{a_{10}}^{2} - 96\,\mathit{t_0}
^{6}\,\mathit{a_{10}} + 360\,\mathit{t_0}^{4}\,\mathit{a_{10}} - 474\,
\mathit{a_{10}}\,\mathit{t_0}^{2} \\
\mbox{} + 216\,\mathit{a_{10}} + 48\,\mathit{t_0}^{8} - 240\,\mathit{
t_0}^{6} + 469\,\mathit{t_0}^{4} - 420\,\mathit{t_0}^{2} + 144) =0
\end{multline*}
Similarly $\cB'$ is described in $\cB$ by the equation
\begin{multline*}
(27 - 45\,\mathit{a_{20}} + 9\,\mathit{a_{20}}^{2} + 9\,\mathit{a_{20}}^{3
} + 19\,\mathit{t_1}^{2} + 18\,\mathit{t_1}^{2}\,\mathit{a_{20}} + 27
\,\mathit{t_1}^{2}\,\mathit{a_{20}}^{2} + 9\,\mathit{t_1}^{4} + 27\,
\mathit{t_1}^{4}\,\mathit{a_{20}} \\
\mbox{} + 9\,\mathit{t_1}^{6})(\mathit{t_1}^{4}\,\mathit{a_{20}}^{2}
 + 3\,\mathit{a_{20}}^{2} - 6\,\mathit{a_{20}} - 6\,\mathit{t_1}^{2}\,
\mathit{a_{20}} + 2\,\mathit{t_1}^{6}\,\mathit{a_{20}} + 2\,\mathit{t_1}
^{4}\,\mathit{a_{20}} + 6\,\mathit{t_1}^{2} + \mathit{t_1}^{8} \\
\mbox{} + 4\,\mathit{t_1}^{4} + 3 + 2\,\mathit{t_1}^{6}) =0
\end{multline*}
One can check that the generic sextics in $\cA'$ is contained in
$\cM(\Si_{13})$,
putting explicit values to parameters.
It is obvious that $\xi_1(\cA')\subset \cB'$ and $\xi_2(\cB')\subset
\cA'$.
Thus the assertion follows.

Next we consider
 the equality  $\psi_{red}(\cM(\Si_{29}))=\psi_{red}(\cM(\Si_{60}))$
with $\Si_{29}=[[4A_2,A_5],[A_5]]$ and $\Si_{60}=[[2A_2,2A_5],[2A_2]]$.
Consider the slice  $\cA'', \cB''$ of 
$\cM(\Si_{29})$ and  $\cM(\Si_{60})$ 
subject to the slice condition $(\star_1+\star_2)$ or $(\star_3)$.
Then we have the canonical inclusions
$ \cA''\subset \cA,\quad \cB''\subset \cB$.
The slices $\cA'',\cB''$ are at the ``boundary'' of $\cA,\cB$
respectively.
For example, consider $(f_2,f_3)\in \cA$. Then if the sextics
$f_2^3+f_3^2=0$ has one inner
$A_5$, it must be on $x$-axis. Thus this is the case
if and only if the resultant $S(y)$ of $f_2(x,y)$ and $f_3(x,y)$
in $x$, which is an even polynomial in $y$, has $y=0$
as a solution. 
%This is the result of the fact that $f_2,f_3$ are
%even in $y$.
This condition  is described as
%$(f_2,f_3)\in cA''$ if and only if 
\begin{multline*}
 - 96\,\mathit{t_0}^{6}\,\mathit{a_{10}} - 180\,\mathit{a_{10}}^{2}\,
\mathit{t_0}^{2} - 24\,\mathit{a_{10}}^{3}\,\mathit{t_0}^{2} + 72\,
\mathit{t_0}^{4}\,\mathit{a_{10}}^{2} - 240\,\mathit{t_0}^{6} + 469\,
\mathit{t_0}^{4} + 360\,\mathit{a_{10}}\,\mathit{t_0}^{4} \\
\mbox{} + 3\,\mathit{a_{10}}^{4} - 420\,\mathit{t_0}^{2} + 216\,
\mathit{a_{10}} + 48\,\mathit{t_0}^{8} + 144 + 30\,\mathit{a_{10}}^{3}
 + 120\,\mathit{a_{10}}^{2} - 474\,\mathit{a_{10}}\,\mathit{t_0}^{2} =0
\end{multline*}
Now we consider $\cB''$. Take  $(g_2,g_3)\in \cB$ and let $S(y)$
be the resultant of $g_2$ and $g_3$ in $x$-variable. As it has an inner
$A_5$ at $O$, $S(y)$ is divisible by $y^2$. 
The condition that the sextics $g_2^3+g_3^2=0$ has two inner $A_5$
is equivalent to $S(y)$ is  divisible by $y^4$. Thus
the slice $\cB''$ consist of $(g_2,g_3)\in \cB$
which satisfy
\begin{multline*}
\mathit{t_1}^{8} + 2\,\mathit{a_{20}}\,\mathit{t_1}^{6} + 2\,\mathit{
t_1}^{6} + \mathit{a_{20}}^{2}\,\mathit{t_1}^{4} + 4\,\mathit{t_1}^{4}
 + 2\,\mathit{a_{20}}\,\mathit{t_1}^{4} - 6\,\mathit{t_1}^{2}\,
\mathit{a_{20}} + 6\,\mathit{t_1}^{2} + 3\,\mathit{a_{20}}^{2} - 6\,
\mathit{a_{20}} \\
\mbox{} + 3 =0
\end{multline*}
Then  after checking that a generic sextics of
 $\cA'',\cB''$
 have the prescribed singularities, the assertion follows from 
$ \xi_1(\cA'')\subset \cB'',\quad \xi_2(\cB'')\subset \cA''$.

\noindent
{\bf III}. We show that  the coincidence of moduli spaces
$ \psi_{red}(\cM(\Si_8))=\psi_{red}(\cM(\Si_{52}))$
where $\Si_8=[[6A_2],[E_6]]$ and $\Si_{52}=[[4A_2,E_6],[2A_2]]$.
First we observe that
$\emdim(\Si_6)=\emdim(\Si_{52})=1$.
In fact, it is easy to see that both moduli spaces are irreducible and
have
transverse slice.
We consider the minimal slices $\cS_{52}$ of  $\cM(\Si_{52})$
(respectively $ \cS_8$  of  $\cM(\Si_8)$) with respect to the slice condition:

($\star$) : two outer (resp. inner) $A_2$ are at $P=(1,1)$ and $Q=(1,-1)$ and 
an inner (resp. outer) $E_6$ is at $O=(0,0)$. The tangent
cones at $P$ and $O$ are given by $y=1$ and $x=0$ respectively.

The normal forms of the slice $\cS_{52}$ and $\cS_8$   can be  given as follows.
\begin{eqnarray*}
&\cS_{52}: \begin{cases}
&\mathrm{f_2}(x, \,y)=y^{2} + ( - 3 - \mathit{t_1}^{2})\,x^{2} + 
2\,x\\
&\mathrm{f_3}(x, \,y)= - {\displaystyle \frac {1}{2}} \,
{\displaystyle \frac { - 6\,y^{2}\,x - 3\,y^{2}\,x\,\mathit{t_1}^{
2} + 6\,y^{2} + 6\,x^{3} + 9\,x^{3}\,\mathit{t_1}^{2} + 2\,x^{3}\,
\mathit{t_1}^{4} - 6\,x^{2}\,\mathit{t_1}^{2} - 6\,x^{2}}{\mathit{
t_1}}}
\end{cases}
\\
&\cS_8:
\begin{cases}
&\mathit{g_2} := y^{2} + (3 - \mathit{t_0}^{2})\,x^{2} + ( - 4 + 2
\,\mathit{t_0}^{2})\,x - \mathit{t_0}^{2}
\\
&\mathit{g_3} := {\displaystyle \frac {1}{2}} (6\,y^{2}\,x - 3\,y
^{2}\,x\,\mathit{t_0}^{2} + 3\,y^{2}\,\mathit{t_0}^{2} + 6\,x^{3}
 - 9\,x^{3}\,\mathit{t_0}^{2} + 2\,x^{3}\,\mathit{t_0}^{4} - 12\,x
^{2} + 21\,x^{2}\,\mathit{t_0}^{2} \\
&\mbox{} - 6\,x^{2}\,\mathit{t_0}^{4} - 12\,x\,\mathit{t_0}^{2} + 6
\,x\,\mathit{t_0}^{4} - 2\,\mathit{t_0}^{4})/\mathit{t_0} 
\end{cases}
\end{eqnarray*}
We can see that $f_2^3+f_3^2=g_2^3+g_3^2$  under the correspondence
 $t_0=2\,I\,\sqrt{3}/t_1$.

\vspace{.3cm}\noindent 
{\bf IV}. We  have already seen
  the coincidence of the images of three moduli spaces
$\cM([[4A_2,A_5],[E_6]])$, $\cM([[4A_2,E_6],[A_5]])$ and
$\cM([[2A_2,A_5,E_6],[2A_2]])$ in the previous section
(Proposition \ref{Triple} ).

\vspace{.3cm}\noindent
{\bf V}. We show that 
$\psi_{red}(\cM(\Si_{123}))=\psi_{red}(\cM(\Si_{t11}))$
where $\Si_{123}=[[2A_2,C_{3,9}^\natural],[A_2]]$ and $\Si_{t11}=[3A_2,C_{3,9}]$.
By  Maple   computation, we can show that
both moduli spaces are irreducible and 
the dimensions of minimal moduli slices are 1.
First, we consider
the minimal moduli slices $\cA$  of  $ \cM (\Si_{123})$
and $\cB$ of  $\cM(\Si_{t11})$ with
three singularities  are specialized as follows.
% Let $O:=(0,0),P:=(0,1), Q:=(1,-1)$. 

$(\star_1)$ for $\cA$:   $C_{3,9}$-singularity is at $(0,0)$ 
with $y=0$ as the tangent line,
$P:=(0,1)$ is an outer  $A_2$-singularity with $x=0$ as the tangent line and 
$Q:=(1,-1)$ is an inner $A_2$-singularity.

$(\star_2)$ for $\cB$:  $C_{3,9}$-singularity is at $(0,0)$ 
with $y=0$ as the tangent line,
$P:=(0,1)$ is an inner $A_2$-singularity with $x=0$ as the tangent line and 
$Q:=(1,-1)$ is an inner $A_2$-singularity.

The normal form of $\cA$ is given by
\begin{multline*}
\mathrm{f_2}(x, \,y)=y^{2} - y - y\,\mathit{t_0}^{2} - y\,x\,
\mathit{t_0}^{2} + {\displaystyle \frac {1}{3}} \,I\,y\,x\,
\mathit{t_0}^{2}\,\sqrt{3} - 2\,x^{2}\,\mathit{t_0}^{2} - 2\,x^{2}
 + {\displaystyle \frac {1}{3}} \,I\,x^{2}\,\mathit{t_0}^{2}\,
\sqrt{3}\\
\mathrm{f_3}(x, \,y)={\displaystyle \frac {1}{4}} (2\,y^{3}\,
\mathit{t_0}^{2} + 3\,I\,x^{3}\,\mathit{t_0}^{2}\,\sqrt{3} - 13\,y
\,x^{2}\,\mathit{t_0}^{2} - 6\,y^{2}\,\mathit{t_0}^{2} - 5\,x^{3}\,
\mathit{t_0}^{2} + I\,y\,x\,\mathit{t_0}^{2}\,\sqrt{3} \\
\mbox{} + I\,y^{2}\,x\,\mathit{t_0}^{2}\,\sqrt{3} + 3\,I\,y\,x^{2}
\,\mathit{t_0}^{2}\,\sqrt{3} - 3\,y\,x\,\mathit{t_0}^{2} - 3\,y^{2}
\,x\,\mathit{t_0}^{2} + 6\,y^{3} - 6\,x^{3} - 6\,y^{2} \\
\mbox{} - I\,y^{2}\,x\,\sqrt{3} - 3\,y\,x + I\,y\,x\,\sqrt{3} + 3
\,y^{2}\,x + 2\,I\,x^{3}\,\sqrt{3} - 12\,y\,x^{2})\mathit{t_0} 
\end{multline*}
We observe that $f_3(x,y)=0$ has a node at $O$ and 
the intersection number $I(C_2,C_3;\, fg
O)=4$. See also \cite{Pho}.
Now we look for another torus decomposition $f(x,y)=g_2(x,y)^3+g_3(x,y)^2$
so that $I(g_2,g_3; O)=3$ and thus the conic $g_2(x,y)=0$ passes through 
three $A_2$ singularities and $C_{3,9}$-singularity. 
By an easy computation, we find $g_2$ is given by
\begin{multline*}
\mathrm{g_2}(x, \,y)= - {\displaystyle \frac {1}{6}} \,
( - 6\,y^{2} - 6\,y^{2}\,\mathit{t_0}^{2} + 6
\,y + 6\,y\,\mathit{t_0}^{2} - 3\,y\,x\,\mathit{t_0}^{2} + I\,y\,x
\,\mathit{t_0}^{2}\,\sqrt{3} \\
+ I\,x^{2}\,\mathit{t_0}^{2}\,\sqrt{3}
 + 9\,x^{2}\,\mathit{t_0}^{2} + 12\,x^{2})/(1 + \mathit{t_0}^{2})
\end{multline*}
To look for a partner cubic form $g_3$, we apply the degeneration method
to the family
$g_t:=f-t\,g_2 \, ^3$.
We can take $t=(1+t_0^2)^3$
and the partner cubic form is given by
%$g_3(x,y)=g_{3,0}s$ where
\begin{multline*}
\mathrm{g_3}(x, \,y)={\displaystyle \frac {1}{8}} (10\,x^{3}\,
\mathit{t_0}^{2} + 4\,y\,x^{2}\,\mathit{t_0}^{2} + 2\,I\,y\,x^{2}\,
\mathit{t_0}^{2}\,\sqrt{3} - 6\,y^{2}\,x\,\mathit{t_0}^{2} + 6\,y\,
x\,\mathit{t_0}^{2} - 3\,y^{3}\,\mathit{t_0}^{2} \\
\mbox{} - I\,\sqrt{3}\,y^{3}\,\mathit{t_0}^{2} + 3\,y^{2}\,
\mathit{t_0}^{2} + I\,y^{2}\,\mathit{t_0}^{2}\,\sqrt{3} + 12\,x^{3}
 + 6\,y\,x^{2} + 2\,I\,y\,x^{2}\,\sqrt{3} - 6\,y^{2}\,x + 6\,y\,x
 \\
\mbox{} - 3\,y^{3} - I\,\sqrt{3}\,y^{3} + 3\,y^{2} + I\,y^{2}\,
\sqrt{3})\sqrt{ - 2 + 2\,I\,\sqrt{3}}\,\mathit{t_0} 
\end{multline*}
where $t_0=f_3(0,1)/f_2(0,1)$.
This gives an isomorphism $\phi: \cA\to \cB$
which completes the proof.

\newpage
\subsection{Configuration table }\label{Table}
\begin{center}Table 1
\end{center}

\begin{tabular}[ht]{|c|c|c|c|c|c|}
\hline
$\iv$& No &$\Sigma$ & $[g,\mu^*,n^*,i(C)]$ &Existence? & $\emdim$ \\ \hline
%$\iv(C)$ is [1,1,1,1,1,1]
[1,1,1,1,1,1]
&nt1 & $ [[6 A_2], [A_1]] $ &[3,13,10,18] &$\to$ t3,nt2 & 6\\  
&nt2 & $ [[6 A_2], [A_2]] $ &[3,14,9,16] & $\to$ t8,nt3 & 5\\  
&nt3 & $ [[6 A_2], [A_3]] $ &[2,15,8,12] & $\to$ t9,nt4 & 4\\  
&nt4 & $ [[6 A_2], [A_4]] $ &[2,16,7,9] & $\to$ t10,nt5 & 3\\  
&nt5 & $ [[6 A_2], [A_5]] $ &[1,17,6,6] & $\to$ t11,nt29& 2\\  
&nt6 & $ [[6 A_2], [B_{3,3}]] $ &[1,16,6,6] & $\to$ nt30,nt7& 3 \\  
&nt7 & $ [[6 A_2], [D_5]] $ &[1,17,5,4] & $\to$ nt31,nt8    & 2 \\  
&nt8 & $ [[6 A_2], [E_6]] $ &[1,18,4,2] &$\to$ nt32         & 1 \\  
&nt9 & $ [[6 A_2], [2 A_1]] $ &[2,14,8,12] & $\to$ nt43,nt10    & 5 \\  
&nt10 & $ [[6 A_2], [A_1, A_2]] $ &[2,15,7,10] & $\to$ nt44,nt11& 4 \\  
&nt11 & $ [[6 A_2], [A_1, A_3]] $ &[1,16,6,6] & $\to$nt45,nt12  & 3 \\  
&nt12 & $ [[6 A_2], [A_1, A_4]] $ &[1,17,5,3] & $\to$ nt46,nt13 & 2 \\  
&nt13 & $ [[6 A_2], [A_1, A_5]] $ &[0,18,4,0] & $\to$ nt47      & 1 \\  
&nt14 & $ [[6 A_2], [A_1, B_{3,3}]] $ &[0,17,4,0] & No  & 2 \\  
&nt15 & $ [[6 A_2], [2 A_2]] $ &[2,16,6,8] &$\to$ nt23    & 3\\  
&nt16 & $ [[6 A_2], [A_2, A_3]] $ &[1,17,5,4] & No        & 2\\  
&nt17 & $ [[6 A_2], [A_2, A_4]] $ &[1,18,4,1] & No        & 1\\  
&nt18 & $ [[6 A_2], [2 A_3]] $ &[0,18,4,0] & No& 1\\  
&nt19 & $ [[6 A_2], [3 A_1]] $ &[1,15,6,6] & $\to$ t15,nt20& 4\\  
&nt20 & $ [[6 A_2], [2 A_1, A_2]] $ &[1,16,5,4] &$\to$ nt69,nt21 & 3\\  
&nt21 & $ [[6 A_2], [2 A_1, A_3]] $ &[0,17,4,0] & $\to$ nt70 & 2\\  
&nt22 & $ [[6 A_2], [A_1, 2 A_2]] $ &[1,17,4,2] & $\to$ nt52,nt23& 2\\  
&nt23 & $ [[6 A_2], [3 A_2]] $ &[1,18,3,0] &Max & 1\\  
&nt24 & $ [[6 A_2], [4 A_1]] $ &[0,16,4,0] &$\to$ nt100 & 3\\  
\hline
\end{tabular}

\newpage
\begin{center}Table 2
\end{center}
%\vspace{.2cm}
\begin{tabular}[ht]{|c|c|c|c|c|c|}
\hline
$\iv$& No &$\Sigma$ & $[g,\mu^*,n^*,i(C)]$ &Existence? & $\emdim$\\ 
\hline
%$\iv(C)$ is [1,1,1,1,2]
[1,1,1,1,2]&nt25 & $ [[4 A_2, A_5], [A_1]] $ &[2,14,10,16] & $\to$ t5,nt26 &5 \\  
&nt26 & $ [[4 A_2, A_5], [A_2]] $ &[2,15,9,14] & $\to$t17, nt27 &4 \\  
&nt27 & $ [[4 A_2, A_5], [A_3]] $ &[1,16,8,10] & $\to$ t18,nt28 &3 \\  
&nt28 & $ [[4 A_2, A_5], [A_4]] $ &[1,17,7,7] & $\to$ nt29      &2 \\  
&nt29 & $ [[4 A_2, A_5], [A_5]] $ &[0,18,6,4] & $\to$ nt32,nt47 &1 \\  
&nt30 & $ [[4 A_2, A_5], [B_{3,3}]] $ &[0,17,6,4] & $\to$ nt31  &2 \\  
&nt31 & $ [[4 A_2, A_5], [D_5]] $ &[0,18,5,2] & $\to$ nt32      &1 \\  
&nt32 & $ [[4 A_2, A_5], [E_6]] $ &[0,19,4,0] & Max             &0 \\  
&nt33 & $ [[4 A_2, A_5], [2 A_1]] $ &[1,15,8,10] & t14,nt34     &4 \\  
&nt34 & $ [[4 A_2, A_5], [A_1, A_2]] $ &[1,16,7,8] &$\to$ nt62,nt35 &3 \\  
&nt35 & $ [[4 A_2, A_5], [A_1, A_3]] $ &[0,17,6,4] & $\to$ nt63,nt36&2 \\  
&nt36 & $ [[4 A_2, A_5], [A_1, A_4]] $ &[0,18,5,1] & $\to$ nt64     &1 \\  
&nt37 & $ [[4 A_2, A_5], [2 A_2]] $ &[1,17,6,6] & $\to$ nt60  &2 \\  
%\hline
%[1,1,1,1,2]
&nt38 & $ [[4 A_2, A_5], [A_2, A_3]] $ &[0,18,5,2] &No &1 \\  
&nt39 & $ [[4 A_2, A_5], [A_2, A_4]] $ &[0,19,4,-1] & No&0 \\  
&nt40 & $ [[4 A_2, A_5], [3 A_1]] $ &[0,16,6,4] & $\to$ nt65,nt41&3 \\  
&nt41 & $ [[4 A_2, A_5], [2 A_1, A_2]] $ &[0,17,5,2] & $\to$nt66,nt42&2 \\  
&nt42 & $ [[4 A_2, A_5], [A_1, 2 A_2]] $ &[0,18,4,0] & $\to$ nt67&1 \\  
%\hline
&nt43 & $ [[4 A_2, E_6], [A_1]] $ &[2,15,8,12] & $\to$t6,nt44&4 \\  
&nt44 & $ [[4 A_2, E_6], [A_2]] $ &[2,16,7,10] & $\to$ t20,nt45&3 \\  
&nt45 & $ [[4 A_2, E_6], [A_3]] $ &[1,17,6,6] & $\to$ t21,nt46&2 \\  
&nt46 & $ [[4 A_2, E_6], [A_4]] $ &[1,18,5,3] & $\to$ nt47&1 \\  
&nt47 & $ [[4 A_2, E_6], [A_5]] $ &[0,19,4,0] & Max&0 \\  
&nt48 & $ [[4 A_2, E_6], [B_{3,3}]] $ &[0,18,4,0] &No &1 \\  
&nt49 & $ [[4 A_2, E_6], [2 A_1]] $ &[1,16,6,6] & $\to$nt68,nt50 &3 \\  
&nt50 & $ [[4 A_2, E_6], [A_1, A_2]] $ &[1,17,5,4] & $\to$ nt69,nt51&2 \\  
&nt51 & $ [[4 A_2, E_6], [A_1, A_3]] $ &[0,18,4,0] & $\to$ nt70&1 \\  
&nt52 & $ [[4 A_2, E_6], [2 A_2]] $ &[1,18,4,2] & $\to$ nt67&1 \\  
&nt53 & $ [[4 A_2, E_6], [3 A_1]] $ &[0,17,4,0] & $\to$ nt71&2 \\  
\hline
\end{tabular}

\begin{center}Table 3
\end{center}
%\vspace{-.2cm}
\begin{tabular}[ht]{|c|c|c|c|c|c|}
\hline
$\iv$& No &$\Sigma$ & $[g,\mu^*,n^*,i(C)]$ &Existence? & $\emdim$\\ \hline
%$\iv(C)$ is [1,1,2,2]
[1,1,2,2]
&nt54 & $ [[2 A_2, 2 A_5], [A_1]] $ &[1,15,10,14] & $\to$ t13&4 \\  
&nt55 & $ [[2 A_2, 2 A_5], [A_2]] $ &[1,16,9,12] & $\to$ nt56&3 \\  
&nt56 & $ [[2 A_2, 2 A_5], [A_3]] $ &[0,17,8,8] & $\to$ nt57 &2 \\  
&nt57 & $ [[2 A_2, 2 A_5], [A_4]] $ &[0,18,7,5] & $\to$ nt64 &1 \\  
&nt58 & $ [[2 A_2, 2 A_5], [2 A_1]] $ &[0,16,8,8] & $\to$ nt96, nt59 &3 \\  
&nt59 & $ [[2 A_2, 2 A_5], [A_1, A_2]] $ &[0,17,7,6] &$\to$ nt97,nt60&2 \\  
%\hline
%[1,1,2,2]
&nt60 & $ [[2 A_2, 2 A_5], [2 A_2]] $ &[0,18,6,4] & $\to$ nt67       &1 \\  
&nt61 & $ [[2 A_2, A_5, E_6], [A_1]] $ &[1,16,8,10] & $\to$ t14,nt62 &3 \\  
&nt62 & $ [[2 A_2, A_5, E_6], [A_2]] $ &[1,17,7,8] & $\to$ nt63      &2 \\  
&nt63 & $ [[2 A_2, A_5, E_6], [A_3]] $ &[0,18,6,4] & $\to$ nt64      &1 \\  
&nt64 & $ [[2 A_2, A_5, E_6], [A_4]] $ &[0,19,5,1] & Max             &0 \\  
&nt65 & $ [[2 A_2, A_5, E_6], [2 A_1]] $ &[0,17,6,4] & $\to$nt98,nt66&2 \\  
&nt66 & $ [[2 A_2, A_5, E_6], [A_1, A_2]] $ &[0,18,5,2] & $\to$ nt67 &1 \\  
&nt67 & $ [[2 A_2, A_5, E_6], [2 A_2]] $ &[0,19,4,0] & Max           &0 \\  
&nt68 & $ [[2 A_2, 2 E_6], [A_1]] $ &[1,17,6,6] & $\to$t15,nt69      &2 \\  
&nt69 & $ [[2 A_2, 2 E_6], [A_2]] $ &[1,18,5,4] & $\to$ nt70         &1 \\  
&nt70 & $ [[2 A_2, 2 E_6], [A_3]] $ &[0,19,4,0] &  Max               &0 \\  
&nt71 & $ [[2 A_2, 2 E_6], [2 A_1]] $ &[0,18,4,0] & $\to$ nt100      &1 \\  
\hline
%$\iv$& No &$\Sigma$ & $[g,\mu^*,n^*,i(C)]$ &Existence? & \\ \hline
%$\iv(C)$ is [1,1,1,3]
[1,1,1,3]&nt72 & $ [[3 A_2, A_8], [A_1]] $ &[2,15,10,15] &$\to$ t19,nt73 &4 \\  
&nt73 & $ [[3 A_2, A_8], [A_2]] $ &[2,16,9,13] & $\to$ t28,nt74  &3 \\  
&nt74 & $ [[3 A_2, A_8], [A_3]] $ &[1,17,8,9] & $\to$ t29,nt75   &2 \\  
&nt75 & $ [[3 A_2, A_8], [A_4]] $ &[1,18,7,6] & $\to $ nt78      &1 \\  
&nt76 & $ [[3 A_2, A_8], [A_5]] $ &[0,19,6,3] &No                &0 \\  
&nt77 & $ [[3 A_2, A_8], [B_{3,3}]] $ &[0,18,6,3] & $\to$ nt78   &1 \\  
&nt78 & $ [[3 A_2, A_8], [D_5]] $ &[0,19,5,1] & Max              &0 \\  
&nt79 & $ [[3 A_2, A_8], [E_6]] $ &[0,20,4,-1] & No              &-1 \\  
&nt80 & $ [[3 A_2, A_8], [2 A_1]] $ &[1,16,8,9] & $\to$ nt108,nt81&3 \\  
&nt81 & $ [[3 A_2, A_8], [A_1, A_2]] $ &[1,17,7,7] & $\to$ nt82   &2 \\  
&nt82 & $ [[3 A_2, A_8], [A_1, A_3]] $ &[0,18,6,3] & $\to$ nt83   &1 \\  
&nt83 & $ [[3 A_2, A_8], [A_1, A_4]] $ &[0,19,5,0] &Max           &0 \\  
%\hline
%[1,1,1,3]
&nt84 & $ [[3 A_2, A_8], [2 A_2]] $ &[1,18,6,5] &No               &1 \\  
&nt85 & $ [[3 A_2, A_8], [A_2, A_3]] $ &[0,19,5,1] &No            &0 \\  
&nt86 & $ [[3 A_2, A_8], [A_2, A_4]] $ &[0,20,4,-2] & No          &-1 \\  
&nt87 & $ [[3 A_2, A_8], [3 A_1]] $ &[0,17,6,3] & $\to$ nt112,nt88&2 \\  
&nt88 & $ [[3 A_2, A_8], [2 A_1, A_2]] $ &[0,18,5,1] & $\to$ nt113&1 \\  
&nt89 & $ [[3 A_2, A_8], [A_1, 2 A_2]] $ &[0,19,4,-1] &No         &0 \\  
%\hline
&nt90 & $ [[3 A_2, B_{3,6}], [A_1]] $ &[0,17,7,6] & $\to$ t20,nt91&3 \\  
&nt91 & $ [[3 A_2, B_{3,6}], [A_2]] $ &[0,18,6,4] & Max           &2 \\  
%\hline
&nt92 & $ [[3 A_2, C_{3,7}], [A_1]] $ &[0,18,6,3] & $\to$ t21     &2 \\  
&nt93 & $ [[3 A_2, C_{3,7}], [A_2]] $ &[0,19,5,1] & No            &1 \\  
\hline
\end{tabular}

\newpage
\begin{center}Table 4
\end{center}\vspace{.2cm}
\begin{tabular}[ht]{|c|c|c|c|c|c|}
\hline
$\iv$& No &$\Sigma$ & $[g,\mu^*,n^*,i(C)]$ &Existence? &$\emdim$ \\ \hline
%
%$\iv(C)$ is [2,2,2]
[2,2,2]&nt94 & $ [[3 A_5], [A_1]] $ &[0,16,10,12] & $\to$ nt95,nt100& 3\\  
&nt95 & $ [[3 A_5], [A_2]] $ &[0,17,9,10] & $\to$ nt97              & 2\\  
&nt96 & $ [[2 A_5, E_6], [A_1]] $ &[0,17,8,8] &$\to$ nt97,nt98      & 2\\  
&nt97 & $ [[2 A_5, E_6], [A_2]] $ &[0,18,7,6] & $\to$ nt99          & 1\\  
&nt98 & $ [[A_5, 2 E_6], [A_1]] $ &[0,18,6,4] & $\to$ nt99,nt100    & 1\\  
&nt99 & $ [[A_5, 2 E_6], [A_2]] $ &[0,19,5,2] &Max                  & 0\\  
&nt100 & $ [[3 E_6], [A_1]] $ &[0,19,4,0] &Max                      & 0\\  
\hline

%$\iv(C)$ is [1,2,3]
[1,2,3]
&nt101 & $ [[A_2, A_5, A_8], [A_1]] $ &[1,16,10,13] & $\to$ nt102 &3 \\  
&nt102 & $ [[A_2, A_5, A_8], [A_2]] $ &[1,17,9,11] & $\to$ nt103  &2 \\  
&nt103 & $ [[A_2, A_5, A_8], [A_3]] $ &[0,18,8,7] & $\to$ nt104   &1 \\  
&nt104 & $ [[A_2, A_5, A_8], [A_4]] $ &[0,19,7,4] &Max            &0 \\  
&nt105 & $ [[A_2, A_5, A_8], [2 A_1]] $ &[0,17,8,7] & $\to$ nt106 &2 \\  
&nt106 & $ [[A_2, A_5, A_8], [A_1, A_2]] $ &[0,18,7,5] & $\to$ nt113&1 \\  
&nt107 & $ [[A_2, A_5, A_8], [2 A_2]] $ &[0,19,6,3] &No             &0 \\  
&nt108 & $ [[A_2, E_6, A_8], [A_1]] $ &[1,17,8,9] &$\to$ nt109   &2 \\  
&nt109 & $ [[A_2, E_6, A_8], [A_2]] $ &[1,18,7,7] & $\to$ nt110  &1 \\  
%\hline
%[1,2,3]
&nt110 & $ [[A_2, E_6, A_8], [A_3]] $ &[0,19,6,3] &Max           &0 \\  
&nt111 & $ [[A_2, E_6, A_8], [A_4]] $ &[0,20,5,0] & No           &-1 \\  
&nt112 & $ [[A_2, E_6, A_8], [2 A_1]] $ &[0,18,6,3] & $\to$ nt113&1 \\  
&nt113 & $ [[A_2, E_6, A_8], [A_1, A_2]] $ &[0,19,5,1] & Max     &0 \\  
&nt114 & $ [[A_2, E_6, A_8], [2 A_2]] $ &[0,20,4,-1] & No        &-1 \\  
\hline
%\end{tabular}
%\begin{tabular}[ht]{|c|c|c|c|c|}
%\hline
%$\iv$& No &$\Sigma$ & $[g,\mu^*,n^*,i(C)]$ & Existence?& \\ \hline
%
[1,1,4]
&nt115 & $ [[2 A_2, A_{11}], [A_1]] $ &[1,16,10,14] &$\to$ t33, nt116 &3 \\  
&nt116 & $ [[2 A_2, A_{11}], [A_2]] $ &[1,17,9,12] & $\to$ nt117      &2 \\  
&nt117 & $ [[2 A_2, A_{11}], [A_3]] $ &[0,18,8,8] & $\to$ nt118       &1 \\  
&nt118 & $ [[2 A_2, A_{11}], [A_4]] $ &[0,19,7,5] &Max                &0 \\  
&nt119 & $ [[2 A_2, A_{11}], [2 A_1]] $ &[0,17,8,8] & $\to$ nt135, nt120&2 \\  
&nt120 & $ [[2 A_2, A_{11}], [A_1, A_2]] $ &[0,18,7,6] & $\to$ nt136    &1 \\  
&nt121 & $ [[2 A_2, A_{11}], [2 A_2]] $ &[0,19,6,4] & No                &0 \\  
%\hline
&nt122 & $ [[2 A_2, C_{3,9}^\natural], [A_1]] $ &[0,18,7,5] & $\to$ t34,nt123 &2 \\  
&nt123 & $ [[2 A_2, C_{3,9}^\natural], [A_2]] $ &[0,19,6,3] & Max             &1 \\  
%\hline
&nt124 & $ [[2 A_2, B_{3,8}], [A_1]] $ &[0,19,6,2] & $\to$ t35     &1 \\  
&nt125 & $ [[2 A_2, B_{3,8}], [A_2]] $ &[0,20,5,0] & No            &0 \\  
\hline
\end{tabular}

\newpage
\begin{center}Table 5
\end{center}\vspace{.2cm}
\begin{tabular}[ht]{|c|c|c|c|c|c|}
\hline
$\iv$& No &$\Sigma$ & $[g,\mu^*,n^*,i(C)]$ &Existence? & $\emdim$\\ \hline

%$\iv(C)$ is [3,3]
[3,3]
&nt126 & $ [[2 A_8], [A_1]] $ &[1,17,10,12] & $\to$ nt127&2 \\  
&nt127 & $ [[2 A_8], [A_2]] $ &[1,18,9,10] & $\to$ nt128 &1 \\  
&nt128 & $ [[2 A_8], [A_3]] $ &[0,19,8,6] &Max           &0 \\  
&nt129 & $ [[2 A_8], [A_4]] $ &[0,20,7,3] & No           &-1 \\  
&nt130 & $ [[2 A_8], [2 A_1]] $ &[0,18,8,6] &$\to$ nt128 &1 \\  
&nt131 & $ [[2 A_8], [A_1, A_2]] $ &[0,19,7,4] & No      &0 \\  
&nt132 & $ [[2 A_8], [2 A_2]] $ &[0,20,6,2] & No         &-1 \\  
\hline

%$\iv(C)$ is [2,4]
%\medskip
[2,4]
&nt133 & $ [[A_5, A_{11}], [A_1]] $ &[0,17,10,12] & $\to$ nt134 &2 \\  
&nt134 & $ [[A_5, A_{11}], [A_2]] $ &[0,18,9,10] & $\to$ nt136  &1 \\  
%\hline
&nt135 & $ [[E_6, A_{11}], [A_1]] $ &[0,18,8,8] & $\to$ nt136   &1 \\  
&nt136 & $ [[E_6, A_{11}], [A_2]] $ &[0,19,7,6] & Max           &0 \\  
\hline

%$\iv(C)$ is [1,5]
%\medskip
%
%\begin{tabular}[ht]{|c|c|c|c|}
%\hline
[1,5]
&nt137 & $ [[A_2, A_{14}], [A_1]] $ &[1,17,10,13] & $\to$ nt138 &2 \\  
&nt138 & $ [[A_2, A_{14}], [A_2]] $ &[1,18,9,11] &$\to$ nt139   &1 \\  
&nt139 & $ [[A_2, A_{14}], [A_3]] $ &[0,19,8,7] &Max            &0 \\  
&nt140 & $ [[A_2, A_{14}], [A_4]] $ &[0,20,7,4] & No            &-1 \\  
&nt141 & $ [[A_2, A_{14}], [2 A_1]] $ &[0,18,8,7] &$\to$ nt142   &1 \\  
&nt142 & $ [[A_2, A_{14}], [A_1, A_2]] $ &[0,19,7,5] &Max        &0 \\  
&nt143 & $ [[A_2, A_{14}], [2 A_2]] $ &[0,20,6,3] & No           &-1 \\ 
\hline

[6]&nt144 & $ [[A_{17}], [A_1]] $ &[0,18,10,12] & $\to$ nt145&1 \\  
&nt145 & $ [[A_{17}], [A_2]] $ &[0,19,9,10] &Max             &0 \\  
\hline
%\end{tabular}
%\bigskip
%

\end{tabular}

\bibliographystyle{abbrv}
%\bibliography{/home/oka/paper/okabib}
%\bibliography{NT-all.bbl}
%\bibliography
%{
%\def\cprime{$'$} \def\cprime{$'$}

%}
\end{document}